\documentclass[12pt, a4paper]{article}
\usepackage{setspace}
\usepackage{geometry}
\usepackage[OT1]{fontenc}
\usepackage{amsthm,amsmath}
\usepackage[colorlinks,citecolor=blue,urlcolor=blue]{hyperref}
\usepackage{amssymb,amsmath,bm,mathrsfs,makeidx,amsfonts,graphicx,amsthm,multirow}
\usepackage{booktabs}

\newcommand\X{\mathbf x}

\def\Im{{\rm\,Im} }
\newcommand\Y{\mathbf y}
\newcommand\bY{\mathbf Y}
\newcommand\bR{{\mathbb R}}
\newcommand\bZ{{\mathbf z}}
\newcommand\Z{{\mathbf Z}}
\newcommand\z{{\mathbf z}}
\newcommand\ee{{\mathbf e}}
\newcommand\y{{\mathbf y}}
\newcommand\w{{\mathbf w}}
\newcommand\W{{\mathbf W}}
\newcommand\wh{\widehat }
\newcommand\cov{\mathrm{Cov}}

\newcommand\pto{\stackrel p\rightarrow}

\newcommand\e{{\mathbb E}}
\newcommand\p{{\mathbb P}}
\newcommand\var{{\rm Var}}
\newcommand\tr{{\mathrm{tr}}}
\newcommand\cF{{\mathcal F}}

\newtheorem{theorem}{Theorem}[section]%
\newtheorem{corollary}[theorem]{Corollary}%
\newtheorem{lemma}[theorem]{Lemma}%
\newtheorem{proposition}[theorem]{Proposition}%

\begin{document}

\begin{center}
\Large The universality principle for spectral distributions of sample covariance matrices.  
\end{center}
\begin{center}
\large Pavel~Yaskov\footnote{Steklov Mathematical Institute, Russia\\
 e-mail: yaskov@mi.ras.ru\\Supported
    by RNF grant 14-21-00162 from the Russian Scientific Fund.}
 \end{center}

\begin{abstract}We derive the universality principle for empirical spectral distributions of  sample covariance matrices  and their Stieltjes transforms. This principle states the following. Suppose  quadratic forms of random vectors $\y_p$ in $\bR^p$ satisfy a weak law of large numbers and the sample size  grows at the same rate as $p$. Then the limiting spectral distribution  of corresponding sample covariance matrices is the same as in the case with conditionally Gaussian  $\y_p$.
This result is generalized for $m$-dependent martingale difference sequences and $m$-dependent linear processes.
\end{abstract}

\begin{center}
\end{center}

\section{Introduction}
The random matrix theory plays an important role in modern high-dimensional statistics (e.g., see \cite{PA}). A large-dimensional sample covariance matrix is an object of  primary interest. Many test statistics could be defined by its eigenvalues (e.g., see \cite{BZ} and \cite{K}). Asymptotic behaviour of such statistics depends on the empirical distribution of the eigenvalues. The latter is called the empirical spectral distribution (ESD). 

The  universality principle for ESDs of sample covariance matrices says that the limiting behaviour of ESDs is the same as when a random sample is taken from a Gaussian distribution. In the pioneering paper \cite{MP}, Marchenko and Pastur discovered general conditions implying universality. Namely, if quadratic forms of random vectors $\y_p$ in $\bR^p$ concentrate near their expectations, then ESDs of corresponding sample covariance matrices obeys the universality principle. Bai and Zhou \cite{BZ} gave the first formal proof of this fact (see also the paper of Pajor and Pastur \cite{PP} and the book of Pastur and Shcherbina \cite{PS}). They assumed that entries of $\y_p$ have finite fourth moments. Girko and Gupta \cite{G}  considered the universality principle without the finite fourth-moment assumption (but under a more restrictive  assumption on covariance matrices).  All of these results are particular cases of a general universality principle that is derived in the present paper. 

Recall that a sample covariance matrix in the random matrix theory is usually defined by   
$n^{-1}\bY_{pn}\bY_{pn}^\top,$
where $\bY_{pn}$ is a $p\times n$ random matrix which columns are independent copies of $\y_p$. The average of these columns is not subtracted since it does not affect the limiting spectral distribution (see Chapter 3 in \cite{BS}). 

The present paper contributes to the random matrix theory in two directions. First, it shows when one has the universality principle for ESDs in a very general framework. 
Namely, if a weak law of large numbers for quadratic forms of $\y_p$ holds,  then  the limiting spectral distribution  of $n^{-1}\bY_{pn}\bY_{pn}^\top$ is the same as in the case of conditionally Gaussian  $\y_p$ when $n$  grows at the same rate as $p$. Similar conditions have appeared in the literature in a much stronger form. E.g., see Theorem 1.1 in \cite{BZ}, Theorem 19.1.8 in \cite{PS} and Theorem 6.1 in \cite{G}.  

We generalize these results for random matrices  $\bY_{pn}$ which columns form an $m$-dependent martingale difference sequence or an $m$-dependent linear process. Recently Banna, Merlevede and Peligrad \cite{mp} obtained the universality principle assuming $m$-dependence in rows and columns of $\bY_{pn}$. However, the technique developed in \cite{mp} allows to derive this property only in the case where the limiting spectral distribution of $n^{-1}\bY_{pn}\bY_{pn}^\top$ is completely determined by the covariance structure of $\bY_{pn}$'s entries. However, in general, it is not determined.  

Second, we derive useful moment inequalities for quadratic forms. These inequalities show when a weak law of large numbers holds. The latter allows us to describe explicitly  a wide class of $\y_p$ which sample covariance matrices have universality properties.

The paper is structured as follows. Sections 2 and 3 contain universality laws for ESDs and their Stieltjes transforms. Sections 4 deals with moment inequalities for quadratic forms. The proofs  are given in Section 5 and Appendix.

\section{Universality of ESD: independent observations.}

For each $p\geqslant 1,$ let $\Y_p$ be a random vector in $\bR^p$ and $\Sigma_p$ be a random symmetric positive semi-definite $p\times p$ matrix defined on the same probability space. Consider the following assumptions.\\
{\bf (A1)} $(\Y_p^\top A_p\Y_p-\tr(\Sigma_p A_p))/p\pto 0$ as $p\to\infty$ for all sequences of real symmetric positive semi-definite $p\times p$  matrices $A_p$ with uniformly bounded spectral norms $\|A_p\|$.\\
{\bf (A2)} $\tr(\Sigma_p^2)/p^2\pto 0$ as $p\to\infty.$

Assumption (A1) is a version of the weak law of large numbers for quadratic forms. It is a key assumption.  Obviously, it is satisfied for $\Sigma_p=\y_p\y_p^\top.$ However, (A2)  may not hold in this case. Assumption (A2) guarantees that (A1) holds when,  conditionally on $\Sigma_p$, $\y_p$ has a centred normal distribution. This is shown in the next proposition (for a proof, see Appendix).
\begin{proposition}\label{pz}
Let $\y_p=\Sigma_p^{1/2}\w_p$, where $\Sigma_p^{1/2}$ is the principal square root of $\Sigma_p$ and $\w_p$ is a standard normal vector in $\bR^p$ that is  independent of $\Sigma_p$ for each $p\geqslant1$. Then $(A1)$ holds if and only if $(A2)$ holds. Moreover,
\[\p(|\Y_p^\top A_p\Y_p-\tr(\Sigma_p A_p)|>\varepsilon p)\leqslant \e\min\Big\{\frac{16M^2\tr(\Sigma_p^2)}{(\varepsilon p)^2},1\Big\}\]
 for any $\varepsilon,M>0$ and each complex $p\times p$ matrix $A_p$ with spectral norm $\|A_p\|\leqslant M.$ 
\end{proposition} 

Let $\mathbf Y_{pn}$ be a $p\times n$ matrix which columns are independent copies of $\y_p$ and 
  $\Z_{pn}$ be a $p\times n$ matrix which columns are independent copies of $\z_p=\Sigma_p^{1/2}\mathbf w_p$, where $\w_p$ is given in Proposition \ref{pz}.   The universality principle for ESDs states that ESD of $n^{-1}\mathbf Y_{pn}\mathbf Y_{pn}^\top$ asymptotically behaves in the same manner as ESD of  $n^{-1}\mathbf Z_{pn}\mathbf Z_{pn}^\top$.  Recall that ESD of a $p\times p$ symmetric matrix $A$ is uniquely defined by its Stieltjes transform in the upper-half plane \[s(z)=p^{-1}\tr(A-zI_p)^{-1},\quad z\in \mathbb C^+=\{w\in\mathbb C:\Im(w)>0\},\] 
where $I_p$ is the $p\times p$ identity matrix. The following theorem  deals with universality properties of Stieltjes transforms of $n^{-1}\mathbf Y_{pn}\mathbf Y_{pn}^\top$. 
\begin{theorem}\label{p1}
 If $(A1)$ and $(A2)$ hold, then
\begin{equation}\label{st}
p^{-1}\tr\big(n^{-1}\mathbf Y_{pn}\mathbf Y_{pn}^\top-zI_p\big)^{-1}-p^{-1}\tr\big(n^{-1}\Z_{pn}\Z_{pn}^\top-zI_p\big)^{-1}\to 0\quad\text{a.s.}
\end{equation} for all $z\in \mathbb C^+$ as $n\to\infty$, where $p=p(n)$ is such that $p/n\to y$ for some $y>0$. 
\end{theorem}

\noindent{\bf Remark 1.} Theorem \ref{p1} can be extended to the case where columns of $\bY_{pn}$ are not identically distributed. Namely,  let  $(\y_{pk},\Sigma_{pk})$  be independent over $k=1,\ldots,n,$ $\bY_{pn}$  and  $\Z_{pn}$ be $p\times n$ matrices with columns $\{\y_{pk}\}_{k=1}^{n}$ and $\{\Sigma_{pk}^{1/2}\w_{pk}\}_{k=1}^{n}$  correspondingly, where $\{\w_{pk}\}_{k=1}^{n}$ are  independent standard normal vectors that are also independent of  $\{\Sigma_{pk}\}_{k=1}^{n}$.
The conclusion of  Theorem \ref{p1} holds when the following assumptions hold for given $p=p(n)$ with $p(n)\to \infty$ as $n\to\infty$.\\
{\bf (A3)} For any $M,\varepsilon>0,$ 
\[\frac{1}{n}\sum_{k=1}^n \sup_{A_p}\p(|\Y_{pk}^\top A_{p}\Y_{pk}-\tr(\Sigma_{pk} A_{p})|>\varepsilon p)\to 0,\quad n\to\infty,\] where each supremum is taken over all real symmetric positive semi-definite $p\times p$  matrices $A_p$ with spectral norms $\|A_p\|\leqslant M.$\\
{\bf (A4)} For any $\varepsilon>0,$  \[\frac{1}{n}\sum_{k=1}^n \p(\tr(\Sigma_{pk}^2)>\varepsilon p^2)\to 0,\quad n\to\infty.\]

  In practice, Theorem \ref{p1} is supposed to be used with the next well-known proposition (e.g., see Exercise 2.4.10 in \cite{Tao} for a more general statement). \begin{proposition}\label{STM}
Let $p=p(n)\to\infty$ as $n\to\infty,$ and  $A_n$ is a real random symmetric $p\times p$ matrix for each $n\geqslant 1.$ If, for each $z\in\mathbb C^+,$
\[p^{-1}\tr\big(A_n-zI_p\big)^{-1}\to s(z)\quad{a.s.}\]
as $n\to\infty$ for a deterministic function $s=s(z)$, then \[s(z)=\int_{-\infty}^\infty\frac{F(d\lambda)}{z-\lambda}\] is the Stieltjes transform of some distribution function $F=F(\lambda)$ on $\bR$ and
\[\p\big(\lim_{n\to\infty}F_n(\lambda)=F(\lambda)\text{ for any continuity point $\lambda$ of }F)=1,\]
where $F_n(\lambda)=\sum_{k=1}^pI(\lambda_{kn}\leqslant \lambda)/p$, $\lambda\in\bR,$ for the set $\{\lambda_{kn}\}_{k=1}^p$  of eigenvalues of $A_n$.
\end{proposition}

Due to Proposition \ref{STM}, to prove that ESD of $n^{-1}\mathbf Y_{pn}\mathbf Y_{pn}^\top$ converges vaguely, one should only check that Stieltjes transforms of $n^{-1}\mathbf Z_{pn}\mathbf Z_{pn}^\top$ from Theorem \ref{p1} converge. 

Let us now consider three important particular cases.

{\bf Case 1. $\Sigma_p$, $p\geqslant 1,$ are deterministic matrices.}  
In this case,
\[\mathbf Z_{pn}\mathbf Z_{pn}^\top=\Sigma_p^{1/2}\W_{pn}\W_{pn}^\top\Sigma_{p}^{1/2},\]
where $\W_{pn}$ is a  $p\times n$ random matrix with independent standard normal entries. The set of eigenvalues of \[\Sigma_p^{1/2}\W_{pn}\W_{pn}^\top\Sigma_{p}^{1/2}\]  is the same as that of $\W_{pn}^\top\Sigma_p\W_{pn}$ (when excluding $|p-n|$ zero  eigenvalues). In addition, since
the distribution $\W_{pn}$ does not change when multiplying $\W_{pn}$ by an orthogonal matrix, 
\[\W_{pn}^\top\Sigma_p\W_{pn}\stackrel{d}=\W_{pn}^\top D_p\W_{pn},\]
where $\stackrel{d}=$ is the equality in law and $D_p$ is a diagonal matrix which diagonal entries are eigenvalues of $\Sigma_p$. Hence, the limiting spectral distribution of $n^{-1}\W_{pn}^\top D_p\W_{pn}$ could be derived from Theorem 4.3 in \cite{BS} (or Theorem 7.2.2 in \cite{PS}), when ESD of $D_p$ converges weakly to a probability measure on $\bR_+$. 

 A very important case is $\Sigma_p=D_p=I_p.$ In order to give some precise statements, we recall that the Marchenko-Pastur law  with parameter $y>0$ has a distribution function 
  \[F_y(\lambda)=\max\{1-1/y,0\}I(\lambda\geqslant 0)+I(\lambda\in[a,b])\int_{a}^\lambda \frac{\sqrt{(b-x)(x-a)}}{2\pi xy}\,dx,\]
  where $\lambda\in\bR$, $a=(1-\sqrt{y})^2$ and  $b=(1+\sqrt{y})^2.$
 \begin{theorem}\label{MP}
 Let $(A1)$ holds for $\Sigma_p=I_p$ and $p=p(n)$. Then, with probability one, 
  ESD of $\mathbf Y_{pn}\mathbf Y_{pn}^\top/n$ converges vaguely to the Marchenko-Pastur law with parameter $y>0$ as $n\to\infty$ and  $p/n\to y$.
\end{theorem}  

{\bf Case 2. $\Sigma_p=\xi I_p$, $p\geqslant1$, for a random variable $\xi$.}  In this case,
\[\mathbf Z_{pn}\mathbf Z_{pn}^\top=\W_{pn}T_n\W_{pn}^\top,\]
where $\W_{pn}$ is as above and $T_n$ is a diagonal matrix with diagonal entries that are independent copies of $\xi^2$. By the Glivenko-Cantelli theorem, ESD of $T_n$ converges weakly a.s. to the distribution of $\xi^2.$ Hence, the limiting spectral distribution of $n^{-1}\W_{pn}^\top T_n\W_{pn}$ can be derived from Theorem 4.3 in \cite{BS} (or Theorem 7.2.2 in \cite{PS}).

{\bf Case 3. $\Sigma_p$, $p\geqslant 1,$ are random diagonal matrices.}  In this case, conditionally on $\{\Sigma_{pk}\}_{k=1}^n$, the matrix $\mathbf Z_{pn}$ consists of centred independent normal variables with  different variances in general. There are some particular results that allow to calculate the limiting spectral distribution  in this case. E.g., see Theorem 2 in \cite{GNT} (or Theorem 1 in \cite{BBS}). However, we are not aware of any general result.

To apply Theorem \ref{p1}, we need tools to  verify (A1). This could be done via moment inequalities for quadratic forms given in Section 4.

\section{Universality of ESD: $m$-dependent observations.}

In this section, we generalize results from Section 1 to the $m$-dependent case. 

For each $p,n\geqslant 1$ and $k=1,\ldots,n,$ let $y_{pk}$ be a random vector in $\bR^p$ and 
$\Sigma_{pk}$ be a symmetric random positive semi-definite $p\times p$ matrix. Define \[\cF_{k}^p=\sigma((y_{pj},\Sigma_{pj}),1\leqslant j\leqslant k),\quad 1\leqslant k\leqslant n,\] $\cF_{0}^p$ to be the trivial $\sigma$-algebra and $\e_{k}=\e(\cdot|\cF_k^p)$ for $k=0,1,\ldots,n$.
 
Assume that $\{(\y_{pk},\Sigma_{pk})\}_{k=1}^n$ is an $m$-dependent sequence for some $m=m(n)$ and $p=p(n)$, i.e. $\sigma$-algebras $\cF_{k}^p$ and $\sigma((\y_{pj},\Sigma_{pj}),k+m\leqslant j\leqslant n)$ are independent for each $k=1,\ldots,n-m.$ Introduce the following assumption.\\
{\bf (A5)} For given $p=p(n)$ and $m=m(n)$, as $n\to\infty$,  
\[\sum\e(\|\e_{k}\Sigma_{pl}\|+\|\e_{k}\y_{pl}\y_{pl}^\top \|)\tr\big(\y_{pk}\y_{pk}^\top+\Sigma_{pk})=o(n^3),\] where the sum is taken over all $(k,l)$ with $1\leqslant k<l\leqslant n$ and $l<k+m.$
Suppose also
\[\sum(\|\e \Sigma_{pl}\|+\|\e \y_{pl}\y_{pl}^\top\|)\tr\big(\e\y_{pk}\y_{pk}^\top+\e\Sigma_{pk})=o(n^3),\]
 where the sum is  taken  over all $(k,l)$ with $1\leqslant k,l\leqslant n$ and $m\leqslant |k-l|<2m.$
 \begin{theorem}\label{p2}
Let $(A3)$, $(A4)$ and $(A5)$ hold. If $m=m(n)$ and $p=p(n)$ satisfy $m=o(n/\log n)$ and $p/n\to y>0$ for some $y>0$ as $n\to \infty$, then
\begin{equation*}
p^{-1}\tr\big(n^{-1}\mathbf Y_{pn}\mathbf Y_{pn}^\top-zI_p\big)^{-1}-p^{-1}\tr\big(n^{-1}\Z_{pn}\Z_{pn}^\top-zI_p\big)^{-1}\to 0\quad\text{a.s.}
\end{equation*} for all $z\in \mathbb C^+$ as $n\to\infty$, where $\Z_{pn}$ and $\bY_{pn}$ are defined in Remark $1$. 
\end{theorem}

We now discuss when all assumption of Theorem \ref{p2} hold. The key assumptions are (A3) and (A5). Assumption (A3) could be verified via moment inequalities for quadratic forms given in Section 4.  Simple conditions implying (A5) are given in the next proposition. 
\begin{proposition}\label{A5} Let 
\[\max_{1\leqslant k\leqslant n}(\|\e_{k-1}(\Sigma_{pk})\|+\|\e_{k-1}(\y_{pk}\y_{pk}^\top)\|)=o(n/m)\;\;\text{a.s.}\] 
for some given $p=p(n)$ and $m=m(n)$.  If, in addition,
\[\sum_{k=1}^n \tr(\e\y_{pk}\y_{pk}^\top+\e\Sigma_{pk})=O(n^2),\]
then $(A5)$ hold. 
\end{proposition}
\noindent{\bf Remark 2.} Assumption (A5) implicitly restricts $\{\y_{pk}\}_{k=1}^n$ to be an {\it almost} martingale difference sequence when $p=p(n)$ grows at the same rate as $n$. 
By {\it almost}, we mean that most entries bof $\e_{k-1}\y_{pk}$ are close to zero for most $k=1,\ldots,n.$ Alternatively, one may say that $\|\e_{k-1}\y_{pk}\|^2=o(n)$ with large probability for most $k$. To see this, consider the following setting. Let $\e\y_{pk}\y_{pk}^\top=I_p$ and $\|\e_{k-1}\y_{pk}\|^2\geqslant Cn$ a.s. for some $C>0$ and all $1\leqslant k \leqslant n$ ($n\geqslant 1$), then
\[\|\e_{k-1}\y_{pk}\y_{pk}^\top\|\geqslant \|\e_{k-1}\y_{pk}(\e_{k-1}\y_{pk})^\top\|=\|\e_{k-1}\y_{pk}\|^2\geqslant Cn,\quad 2\leqslant k \leqslant n,\]
and 
\[\sum_{k=2}^{n}\e \|\e_{k-1}\y_{pk}\y_{pk}^\top \|\tr\big(\y_{p(k-1)}\y_{p(k-1)}^\top)\geqslant Cn^3.\]
Clearly (A5) doesn't hold in this case.
\\
{\bf Example 1.} Suppose $m\geqslant 1$ is fixed and $p=p(n)$. Consider independent identically distributed $\bR^p$-valued random vectors $\ee_{pk},$ $-m<k\leqslant n$, which entries are independent copies of a random variable $e$ with $\e e=0$ and $\e e^2=1$. Let $\y_{pk}=A_{pk}\ee_{pk}$ for all $1\leqslant k\leqslant n$, where $A_{pk}$ is a $p\times p$ random matrix measurable w.r.t. $\sigma(\ee_{pl},\,k-m<l<k)$. 
Let also
\[\frac{1}{n}\sum_{k=1}^n \frac{\tr(\e\Sigma_{pk}^2)}{p^2}=o(1)\quad \text{and}\quad 
\frac{1}{n}\sum_{k=1}^n \frac{\tr(\e\Sigma_{pk})}{p}=O(1),\quad n\to\infty,\]
where $\Sigma_{pk}=A_{pk}A_{pk}^\top$, $1\leqslant k\leqslant n$.  The latter yields (A3) and (A4) (see Corollary \ref{t7} below). Note also that $\e\Sigma_{pk}=\e \y_{pk} \y_{pk}^\top$, $1\leqslant k\leqslant n$. If, in addition, \[\max_{1\leqslant k\leqslant n}\|\Sigma_{pk}\|=o(n)\quad\text{a.s.,}\] then (A5) hold (for $p=O(n)$). This can be proven similarly to Proposition \ref{A5}.

We now extend previous results  to the case of linear $m$-independent processes. As in Section 2, for each $p\geqslant 1,$ let $\Y_p$ be a random vector in $\bR^p$ and $\Sigma_p$ be a random symmetric positive semi-definite $p\times p$ matrix defined on the same probability space.  

Let $\mathbf Y_{pn}$ be a $p\times n$ matrix which columns are independent copies of $\y_p$ and 
  $\Z_{pn}$ be a $p\times n$ matrix which columns are independent copies of $\z_p=\Sigma_p^{1/2}\mathbf w_p$, where $\w_p$ is  a standard normal vector independent of $\Sigma_p$.  Suppose $L_n=(l_{kj})_{k,j=1}^n$ is a lower triangular $m$-banded $n\times n$ matrix for each $n\geqslant 1$ with entries $l_{kj}=l_{kj}(n)$ that are equal to zero when $j\leqslant k-m$ or $j>k$. Then
\[\frac{1}{n}\mathbf Y_{pn}L_n^\top L_n\mathbf Y_{pn}^\top=\frac{1}{n}\sum_{k=1}^n 
\Big(\sum_{j:\, k-m<j\leqslant k}l_{kj}\y_{pj}\Big)\Big(\sum_{j:\, k-m<j\leqslant k} l_{kj} \y_{pj}\Big)^\top\]
with $j\in\{1,\ldots,n\}.$
  Clearly, the sequence
  \[\Big\{\sum_{j:\, k-m<j\leqslant k}l_{kj}\y_{pj}\Big\}_{k=1}^n\]
  is an $m$-dependent linear process. To state the universality principle for these sequences, we need one more assumption.
  \\
  {\bf (A6)} $(\|\e\y_{p}\y_{p}^\top\|+\|\e\Sigma_p\|)\tr(\e\y_{p}\y_{p}^\top+\e\Sigma_{p})=o(p^2)$ as $p\to\infty.$

\begin{theorem}\label{p3}
 If $(A1)$, $(A2)$ and $(A6)$ hold, $m$ is a fixed natural number, and entries of $L_n$ are uniformly bounded over $n$, then
\begin{equation}\label{st}
p^{-1}\tr\big(n^{-1}\mathbf Y_{pn}L_n^\top L_n\mathbf Y_{pn}^\top-zI_p\big)^{-1}-p^{-1}\tr\big(n^{-1}\Z_{pn}L_n^\top L_n\Z_{pn}^\top-zI_p\big)^{-1}\to 0\quad\text{a.s.}
\end{equation} for all $z\in \mathbb C^+$ as $n\to\infty$, where $p=p(n)$ is such that $p/n\to y$ for some $y>0$. 
\end{theorem}
\noindent{\bf Remark 3.} Theorem \ref{p3} can be extended to the case where columns of $\bY_{pn}$ are not identically distributed similarly to Theorem \ref{p1} (see Remark 1). It can be also extended to the case where $m=m(n)$ grows to infinity as $n\to\infty.$

\section{Moment inequalities for quadratic forms}

Let $\{X_k\}_{k=1}^\infty $ be a sequence of  random variables and  $\{\varphi_k\}_{k=1}^\infty$  be a sequence of non-negative numbers such that
\begin{equation}\label{vphi}
|\e X_iX_{j}X_{k}X_{l}|\leqslant \min\{\varphi_{j-i},\varphi_{k-j},\varphi_{l-k}\} \quad\text{for all }i<j<k<l.
\end{equation}
Set $\X_p=(X_1,\ldots,X_p)$ for any $p\geqslant 1$. Consider the following assumption.
\\
{\bf (B1)} $\Phi_0+\Phi_1<\infty$, where $\Phi_0=\sup\{\e X_k^4:k\geqslant 1\}$ and $ \Phi_1=\sum_{k=1}^\infty k\varphi_k.$ 
\begin{theorem}\label{t1}
If $(B1)$ holds, then, for any $a\in\bR^p$ and all real $p\times p$ matrices $A$ with zero diagonal, \[\e|\X_p^\top a|^4\leqslant C(\Phi_0+\Phi_1)\|a\|^4\quad\text{and}\quad \e|\X_p^\top A\X_p|^2\leqslant C(\Phi_0+\Phi_1)\tr(A A^\top)\] for some universal constant  $C>0$, where $\|a\|=\sqrt{a^\top a}$.
\end{theorem} 

The proof of Theorem \ref{t1} is based on the strategy developed by Gaposhkin in \cite{G72}.

\begin{corollary}\label{c1}
If $\y_p=\X_p,$  $p\geqslant 1.$  If $(B1)$ holds, then $(A1)$ and $(A2)$ hold  for diagonal matrices $\Sigma_p$ with diagonal entries $X_1^2,\ldots,X_p^2$.
\end{corollary} 

Let now $\{\phi_k\}_{k=1}^\infty$ be a sequence of non-negative numbers such that 
\begin{equation}\label{phi}
\cov(X_i^2,X_j^2)\leqslant \phi_{j-i} \quad\text{for all }i<j.
\end{equation} 
Introduce the following assumption.\\
{\bf (B2)} $\Phi_0+\Phi_1+\Phi_2<\infty$, where $\Phi_0,$ $\Phi_1$ are given in $(B1)$ and $\Phi_2=\sum_{k=1}^\infty \phi_k$.
\begin{theorem}\label{t2}
If  $(B2)$ holds, then, for all real   $p\times p$ matrices $A$,  
\begin{equation*}
\e|\X_p^\top A\X_p -\tr(\Sigma_p A)|^2\leqslant C(\Phi_0+\Phi_1+\Phi_2)\tr(A A^\top)
\end{equation*}
 for a universal constant $C>0$, where $\Sigma_p$ is a diagonal matrix with diagonal entries $\e X_1^2,\ldots,\e X_p^2$.
\end{theorem} 

\noindent{\bf Remark 4.} Lemma 2.5 in \cite{PP} and Lemma 19.1.4 in \cite{PS} contain some related estimates for isotropic $\X_p$ with a log-concave distribution.

 Since $\tr(A_pA_p^\top)\leqslant p\|A_p\|^2$ for all real $p\times p$ matrices $A_p,$ we have the following version of the law of large numbers under assumptions of Theorem \ref{t2}.

\begin{corollary}\label{c3}
Let $\y_p=\X_p,$  $p\geqslant 1.$ If $(B2)$ holds,  then  $(A1)$ and $(A2)$ hold for the diagonal matrix $\Sigma_p$ with   diagonal entries $\e X_1^2,\ldots,\e X_p^2$.
\end{corollary}

The above theorems assume that all $X_k$ have finite forth moments. It can be a restrictive assumption in some applications. Assuming a martingale-type dependence in $X_k$, one can relax this assumption. \\
\noindent{\bf (B3)} Let $\{X_k\}_{k=1}^\infty$ is a martingale difference sequence  w.r.t. its own filtration and  there is $M>0$ such that $\e X_1^2\leqslant M$ and $\e[X_k^2|X_1,\ldots,X_{k-1}]\leqslant M$  a.s. for all $k\geqslant 2$.
 \begin{theorem}\label{t5}
Let $(B3)$ holds. Then, for all real   $p\times p$ matrices $A$ with zero diagonal,  \begin{equation}\label{lbb}
\e|\X_p^\top A\X_p|^2\leqslant 2M^2\,\tr(AA^\top).\end{equation} 
\end{theorem} 

Let us consider a stronger alternative to $(B3).$\\
\noindent{\bf (B4)} Let $\{X_k\}_{k=1}^\infty$, $\{X_k^2-1\}_{k=1}^\infty$ are martingale difference sequences w.r.t. their own filtrations.

\begin{theorem}\label{t6}
If $(B4)$ holds, then \begin{equation}\label{lb}
\e|\X_p^\top A\X_p -\tr(A)|\leqslant Cb\sqrt{\tr(AA^\top)}+C\sum_{k=1}^p |a_{kk}|\e |X_k^2-1|I(|X_k^2-1|>b^2)\end{equation} for all real  $p\times p$ matrices $A=(a_{ij})_{i,j=1}^p$, each $b>1$ and a universal constant $C>0$.
\end{theorem} 

\noindent{\bf Remark 5.} Using arguments similar to those in the proof of Theorem \ref{t6}, one can derive similar inequalities for  $m$-dependent sequences $\{X_{k}\}_{k=1}^\infty$. The universality principle for such sequences is proved in \cite{HP}.

\begin{corollary}\label{c1}
If $\y_p=\X_p,$  $p\geqslant 1,$ and 
\begin{equation}\label{LC}
\frac{1}{p}\sum_{k=1}^p\e X_k^2I(|X_k|>\varepsilon \sqrt{p})\to 0\quad \text{for all }\varepsilon>0.
\end{equation}
  If $(B4)$ holds, then $(A1)$ and $(A2)$ hold  for diagonal matrices $\Sigma_p=I_p$, $p\geqslant 1$.
\end{corollary} 
Corollary \ref{c1} can be derived from Theorem \ref{t6} by taking $b=\varepsilon \sqrt{p}$, tending  $p\to\infty$ and then tending $\varepsilon\to 0.$ 

Let us now give  examples of $\{X_k\}_{k=1}^\infty$ that satisfy the above assumptions. \\
\noindent {\bf Example 2.} If $\{X_k\}_{k=1}^\infty$ are independent random variables with uniformly bounded forth moments and $\e X_k=0$, then \eqref{vphi} and \eqref{phi} hold for $\varphi_k=\phi_k=0,$ $k\geqslant 1.$\\
\noindent {\bf Example 3.} If $\{X_k\}_{k=1}^\infty$ is a martingale difference sequence with finite  forth moments, then  \eqref{vphi} holds for $\varphi_k=0,$ $k\geqslant 1.$\\
\noindent {\bf Example 4.} Let $\{X_k\}_{k=1}^\infty $ be centred, orthonormal, strongly mixing random variables  with mixing coefficients $(\alpha_k)_{k=1}^\infty$. If these variables have uniformly bounded moments of order $4\delta$ for some $\delta>1$, then \eqref{vphi} and \eqref{phi} hold with \[\varphi_k=\phi_k=C\alpha_k^{(\delta-1)/\delta}\] for large enough $C>0$ (see Corollary A.2 in \cite{HH}). One can give similar bounds for many other weakly dependent sequences.\\
\noindent {\bf Example 5.} If $\{X_k\}_{k=1}^\infty$ are independent identically distributed random variables with $\e X_k=0$ and $\e X_k^2=1$, then  \eqref{vphi} and \eqref{phi} hold for $\varphi_k=\phi_k=0,$ $k\geqslant 1,$ as well as Lindeberg's condition \eqref{LC} hold. 

We  now assume that $\{X_{k}\}_{k=1}^\infty$ is a sequence of orthonormal variables.  Denote by  $\mathcal X$ the set of all random variables \[Y=\sum_{k=1}^\infty c_{k}X_k\quad\text{a.s.}\] for some $\{c_{k}\}_{k=1}^\infty$ with $\sum_{k=1}^\infty c_{k}^2<\infty$, where the first series converges in mean square. Let also 
$\mathcal X_p$ be the set of all $\bR^p$-valued random vectors $\y_p$ which entries belong to $\mathcal X$.
\begin{corollary}\label{t4}
Let $\y_p\in\mathcal X_p,$ $\Sigma_p=\e\y_p\y_p^\top$ and  $A_p$ be a real  symmetric positive semi-definite $p\times p$  matrix for $p\geqslant 1$.
If $(B2)$ holds, then \[\e|\y_p^\top A_p\y_p -\tr(\Sigma_pA_p)|^2\leqslant C(\Phi_0+\Phi_1+\Phi_2)\|A_p\|^2\tr( \Sigma_p^2)\] 
 for a universal constant $C>0$. 
Moreover, if $(A2)$ holds, then $(A1)$ holds.
\end{corollary}  

\begin{corollary}\label{t7}
Let $\y_p\in\mathcal X_p,$ $\Sigma_p=\e\y_p\y_p^\top$ and  $A_p$ be a real  symmetric positive semi-definite $p\times p$  matrix for some $p\geqslant 1$.
If $(B4)$ holds, then
\[\e|\y_p^\top A_p\y_p -\tr(\Sigma_pA_p)|\leqslant Cb\|A_p\|\sqrt{\tr(\Sigma_p^2)}+CL(b)\|A_p\|\tr(\Sigma_p)\] 
 for a universal constant $C>0$ and $b>1$,  where  
 \[L(b)=\sup_{k\geqslant 1}\e |X_k^2-1|I(|X_k^2-1|>b^2).\]
 Moreover, if $(A2)$ holds, $\tr(\Sigma_p)=O(p)$ as $p\to\infty$ and $\{X_k^2\}_{k=1}^\infty$ is a uniformly integrable family, then $(A1)$ holds.
\end{corollary}

\noindent {\bf Remark 6.}  Papers \cite{BM}, \cite{PM}, \cite{PSc}  and \cite{Y} contain the universality principle for  some $\y_p\in\mathcal X_p$ when $\{X_k\}_{k=1}^\infty$ are independent identically distributed random variables with $\e X_k=0,$ $\e X_k^2=1$. This could be derived from the general universality principle given in the present paper.\\

\section{Proofs}

{\bf Proof of Theorem \ref{p1}.} Fix $z\in\mathbb C^+.$ First we proceed as in Step 1 of the proof of Theorem 1.1 in \cite{BZ} (see also the proof of (4.5.6) on page 83 in \cite{BS}) to show that  $S_n(z)-\e S_n(z)\to 0$ a.s. as $n\to\infty,$  where
\[S_n(z)=p^{-1}\tr\big(n^{-1}\bY_{pn}\bY_{pn}^\top-zI_p\big)^{-1}.\]  Similar arguments yield that $s_n(z)-\e s_n(z)\to 0$ a.s. as $n\to\infty,$ where \[s_n(z)=p^{-1}\tr\big(n^{-1}\Z_{pn}\Z_{pn}^\top-zI_p\big)^{-1}.\] Hence, we only need to show that  $\e S_n(z)-\e s_n(z)\to 0$. 

We will use  Lindeberg's method as in the proof of  Theorem 6.1 in \cite{G}. Let $\Y_{p1},\ldots,\Y_{pn}$ and $\bZ_{p1},\ldots,\bZ_{pn}$ be
columns of $\bY_{pn}$ and  $\Z_{pn}$ correspondingly. Assume also w.l.o.g. that  $\{(\y_{pk},\Sigma_{pk})\}_{k=1}^n$ are independent copies of $(\y_{p},\Sigma_{p})$ and $\z_{pk}=\Sigma_{pk}^{1/2}\w_{pk}$ for all $1\leqslant k\leqslant n$, where $\{\w_{pk}\}_{k=1}^n$ are independent standard normal vectors in $\bR^p$ that are also independent of  $\{(\y_{pk},\Sigma_{pk})\}_{k=1}^n$.

 Recall that 
\[\bY_{pn}\bY_{pn}^\top=\sum_{k=1}^n\Y_{pk}\Y_{pk}^\top\quad\text{and}\quad
\Z_{pn}\Z_{pn}^\top=\sum_{k=1}^n\bZ_{pk}\bZ_{pk}^\top.\]
Using this representation, we derive that 
\begin{align*}
S_n(z)&=\frac{1}{p}\tr\Big(n^{-1}\sum_{k=1}^n\y_{pk}\y_{pk}^\top-z I_p\Big)^{-1},\\
 s_n(z)&=\frac{1}{p}\tr\Big(n^{-1}\sum_{k=1}^n\z_{pk}\z_{pk}^\top-z  I_p\Big)^{-1}
\end{align*}
and $|S_n(z)-s_n(z)|\leqslant \sum_{k=1}^n|I_{kn}|/p,$ where 
\begin{align*}
I_{kn} =\tr\Big(C_{kn}+\frac{\y_{pk}\y_{pk}^\top }n-z I_p\Big)^{-1}-\tr\Big(C_{kn}+ \frac{\z_{pk}\z_{pk}^\top }{n}-z I_p\Big)^{-1}
\end{align*}
for
$C_{1n}=\sum_{i=2}^n \z_{pi}\z_{pi}^\top /n,$ $C_{nn}=\sum_{i=1}^{n-1} \y_{pi}\y_{pi}^\top /n$ and
  \[C_{kn}=\frac{1}{n}\sum_{i=1}^{k-1} \y_{pi}\y_{pi}^\top+\frac{1}{n}\sum_{i=k+1}^n \z_{pi}\z_{pi}^\top ,\quad 1<k<n.\]
By  the Sherman-Morrison-Woodbury formula,
\[\tr(C+ww^\top-zI_p)^{-1}-\tr(C-zI_p)^{-1}=-\frac{w^\top (C-zI_p)^{-2}w}{1+w^\top (C-zI_p)^{-1}w}\] for any real symmetric  $p\times p$ matrix $C$ and $w\in\bR^p$. Hence, adding and subtracting $\tr(C_{kn}-z I_p)^{-1}$ to $I_{kn}$ yield
\[I_{kn} =-\frac{ \y_{pk}^\top A_{kn}^{2}\y_{pk}/n}{1+ \y_{pk}^\top A_{kn}\y_{pk}/n}+\frac{ \z_{pk}^\top A_{kn}^2\z_{pk}/n}{1+\z_{pk}^\top A_{kn}\z_{pk}/n},\] where we set $A_{kn}=A_{kn}(z)=(C_{kn}-zI_p)^{-1}$, $1\leqslant k\leqslant n.$

Let us now show that \[\frac1p\sum_{k=1}^n\e|I_{kn}|\to 0.\]
The latter implies that $\e S_n(z)-\e s_n(z)\to 0$.  To estimate $I_{kn}$ we need the following lemma which proof is given in Appendix.

\begin{lemma} \label{l1} Let $w\in\mathbb \bR^p,$ $C$ be a $p\times p$ real symmetric matrix and  $z\in\mathbb C$ with $\Im(z)>0.$ Then 
\[\frac{|w^\top (C-zI_p)^{-2} w|}{|1+w^\top (C-zI_p)^{-1} w|}\leqslant \frac{1}{\Im(z)}.\]
\end{lemma}

 Write $z=u+iv$ for $u\in\bR$ and $v=\Im(z)>0$.
By Lemma \ref{l1}, $|I_{kn}|\leqslant 2/v$. Denote $x\wedge y=\min\{x,y\}.$ Using inequalities \begin{equation}\label{minsum}
|x+y|\wedge 1\leqslant(|x|+|y|)\wedge 1 \leqslant |x|\wedge 1+|y|\wedge 1,\quad x,y\in\mathbb C,\end{equation} we derive  that
\[\e|I_{kn}|=\e(|I_{kn}|\wedge(2/v))\leqslant\e(|\Delta_{kn}|\wedge(2/v))+\e(|\widehat \Delta_{kn}|\wedge(2/v))\]
where
\[\Delta_{kn}=\frac{ \y_{pk}^\top A_{kn}^2\y_{pk}/n}{1+ \y_{pk}^\top A_{kn}\y_{pk}/n}-\frac{\tr(\Sigma_{pk} A_{kn}^2)/n}{1+\tr(\Sigma_{pk} A_{kn})/n},\]\[\widehat \Delta_{kn}=\frac{ \z_{pk}^\top A_{kn}^2\z_{pk}/n}{1+ \z_{pk}^\top A_{kn}\z_{pk}/n}-\frac{\tr(\Sigma_{pk} A_{kn}^2)/n}{1+\tr(\Sigma_{pk} A_{kn})/n}.\]
We estimate only  the term $\e(|\Delta_{kn}|\wedge(2/v))$, since $\e(|\widehat \Delta_{kn}|\wedge(2/v))$ can be estimated similarly (via Proposition \ref{pz}).

Fix any $k\in\{1,\ldots,n\}$. For notational simplicity, we will further write \[\y_p,A_n,\Delta_n,C_n,\Sigma_{p}\quad \text{instead of}\quad \y_{pk},A_{kn},\Delta_{kn},C_{kn},\Sigma_{pk}\] and use the following properties:   $C_n$ is a real symmetric positive semi-definite $p\times p$ random matrix, $(\y_p,\Sigma_p)$ is independent of $A_n=(C_n-zI_p)^{-1}$.

\begin{lemma} \label{l0} Let $A=(C-zI_p)^{-1}$, where $C$ is a real symmetric positive semi-definite $p\times p$ matrix and $z\in\mathbb C^+$. Then the spectral norm of $A$ satisfies $\|A\|\leqslant 1/\Im(z)$.
\end{lemma}

\begin{lemma} \label{l2} Let $(A1)$ holds.  Then, for each $\varepsilon,M>0$,
\begin{equation}\label{limsup}
\lim_{p\to\infty}\sup_{A_p}\p(|\y_p^\top A_p \y_p-\tr(\Sigma_pA_p)|>\varepsilon p)=0,\end{equation}
where the supremum is taken over all complex $p\times p$ matrices $A_p$ with $\|A_p\|\leqslant M.$
\end{lemma}

Fix any $\varepsilon>0$. Take
\[D_{n}=\bigcap_{j=1}^2\{|\y_p^\top (A_{n})^j \y_p-\tr(\Sigma_p (A_{n})^j)|\leqslant \varepsilon p\}\]
 and derive that
\[\e(|\Delta_{n}|\wedge(2/v))\leqslant \e(|\Delta_{n}|\wedge(2/v))I(D_{n})+2\p\big(\overline{D}_{n}\big)/v.\]
By the law of iterated mathematical expectations and Lemma \ref{l0},
\[\p\big(\overline{D}_{n}\big)=\e\big[\p\big(\overline{D}_{n}|A_{n}\big)\big]\leqslant 2\sup_{\widehat A_p}\p(|\y_p^\top \widehat A_p \y_p-\tr(\Sigma_p \widehat  A_p)|>\varepsilon  p),\]
where $ \widehat  A_p$ are as in Lemma \ref{l2} with $M=\max\{v^{-1},v^{-2}\}$.
\begin{lemma} \label{l3} Let $z_1,z_2,w_1,w_2\in\mathbb C$. If $|z_1-z_2|\leqslant \gamma,$ $|w_1-w_2|\leqslant \gamma,$
\[\frac{|z_1|}{|1+w_1|}\leqslant M\] and $|1+w_2|\geqslant\delta$ for some $\delta,M>0$ and $\gamma\in(0,\delta/2)$, then 
\[\Big|\frac{z_1}{1+w_1}-\frac{z_2}{1+w_2}\Big|\leqslant C\gamma\]
for some $C=C(\delta,M)>0.$
\end{lemma}

\begin{lemma} \label{l4} Let $z\in\mathbb C^+$, $\Sigma$ and $C$ be real symmetric positive semi-definite  $p\times p$ matrices. Then
\[|1+\tr(\Sigma (C-zI_p)^{-1})|\geqslant \frac{\Im(z)}{|z|}.\]
\end{lemma}

By Lemma \ref{l4}, we get
\begin{equation}\label{1w}
|1+\tr(\Sigma_p A_{n})/n|=|1+\tr((\Sigma_p/n) (C_n-zI_p)^{-1})|\geqslant \delta=\frac{v}{|z|}.
\end{equation}
Take $\gamma=\varepsilon p/n,$
\[(z_1,w_1)=(\y_p^\top A_{n}^2 \y_p, \y_p^\top A_{n} \y_p)/n,\qquad 
(z_2,w_2)=(\tr(\Sigma_p A_{n}^2),\tr(\Sigma_p A_{n}))/n\]
in Lemma \ref{l3}. By Lemma \ref{l1}, 
\[\frac{|z_1|}{|1+w_1|}\leqslant \frac{1}{v}.\]
By \eqref{1w}, 
\[|1+w_2|\geqslant \delta>3\varepsilon y>\frac{2\varepsilon p}n\]
for small enough $\varepsilon>0$ and large enough $p$ (since $p/n\to y>0$).

Using Lemma \ref{l3}, we derive
\[\e(|\Delta_{n}|\wedge (2/v)) I(D_n)\leqslant \e|\Delta_{n}| I(D_n)\leqslant C(\delta,1/v)\,\frac{\varepsilon p}{n}.\]
Combining all above estimates together yields
\[ \e(|\Delta_{kn}|\wedge(2/v))\leqslant C(\delta,1/v)\frac{\varepsilon p}{n}+ \frac 4v\,\sup_{\widehat A_p}\p(|\y_p^\top \widehat A_p \y_p-\tr(\Sigma_p \widehat  A_p)|>\varepsilon  p)
\]
for each $k=1,\ldots,n$ and
\begin{align*} \frac1p\sum_{k=1}^n\e(|\Delta_{kn}|\wedge(2/v))\leqslant C(\delta,1/v)\varepsilon + \frac {4n}{vp}\,\sup_{\widehat A_p}\p(|\y_p^\top \widehat A_p \y_p-\tr(\Sigma_p \widehat  A_p)|>\varepsilon  p).\end{align*}
Taking $\varepsilon$ small enough and then $p$ large enough, we can make the right hand side of the last inequality arbitrarily small by Lemma \ref{l2} (recall also that $n/p\to1/ y>0$). Thus,
\[\lim_{n\to\infty}\frac1p\sum_{k=1}^n\e(|\Delta_{kn}|\wedge(2/v))=0.\] 

Arguing as above with the help of Proposition \ref{pz}, one can prove that 
\[\lim_{n\to\infty}\frac1p\sum_{k=1}^n\e(|\widehat\Delta_{kn}|\wedge(2/v))=0.\] 
This finishes the proof.  Q.e.d.

\noindent {\bf Proof of Proposition \ref{STM}}. Using the Vitali convergence theorem (see Lemma 2.14 on page 37 in \cite{BS}), one can prove that 
\begin{equation}\label{sns}
\p(s_n(z)\to s(z)\quad\text{for all $z\in\mathbb C^+$})=1,
\end{equation}
where $s_n(z)=p^{-1}\tr(A_n-zI_p)^{-1}.$ For a proof, see Step 3 of the proof of Theorem 2.9 on page 37 in \cite{BS}.
Having \eqref{sns}, it is easy to finish the proof by applying Theorem B.9 in \cite{BS}.  Q.e.d.

\noindent {\bf Proof of Theorem \ref{MP}.} The result follows from  Theorem \ref{p1} and Proposition \ref{STM} as well as Theorem 1.1 in \cite{BZ}.   Q.e.d.

\noindent {\bf Proof of Theorem \ref{p2}.}  The proof is along the proof of Theorem \ref{p1}. However, additional arguments are needed.  In what follows, we will use the same notations as in the proof of Theorem \ref{p1}.

 First, to prove that $S_n(z)-\e S_n(z)\to 0$ a.s., we use results of \cite{mp}. Namely, by A.2 in Chapter 9 of \cite{AMO}, non-zero eigenvalues of
\[\mathbf X_{np}=\frac{1}{\sqrt{p+n}}\begin{pmatrix}
O_p& \bY_{pn}\\
\bY_{pn}^\top& O_n
\end{pmatrix}\]
are square roots of non-zero eigenvalues of  $n^{-1} \bY_{pn} \bY_{pn}^\top$ multiplied by $n/(p+n)$ and their negatives, where $O_k$ is the $k\times k$ zero matrix for $k\geqslant 1$. If the second matrix has exactly  $j$ zero eigenvalues, then
\begin{align*}
S_n(z)&=\frac{1}{p}\sum_{k=1}^{p-j}\frac{1}{\lambda_{kn}-z}-\frac{j}{zp}=
\frac{1}{2p\sqrt{z}}\sum_{k=1}^{p-j}\bigg[\frac{1}{\sqrt{\lambda_{kn}}-\sqrt{z}}+\frac{1}{-\sqrt{\lambda_{kn}}-\sqrt{z}}\bigg]-\frac{j}{zp},
\end{align*}
where $\lambda_{kn},$ $k=1,\ldots, p-j$, are non-zero eigenvalues of  $n^{-1} \bY_{pn} \bY_{pn}^\top$ and $\sqrt{z}$ is chosen in a way that $\sqrt{z}\in \mathbb C^+$ for $z\in\mathbb C^+$.
Define 
\[\wh S_n(w)=\frac{1}{p+n}\tr(\mathbf X_{np}-wI_{p+n}\big)^{-1},\quad w\in \mathbb C^+.\]
Then 
\[\wh S_n(w)=\frac{1}{p+n}\sum_{k=1}^{p-j}\bigg[\frac{1}{\sqrt{n\lambda_{kn}/(p+n)}-w}+\frac{1}{-\sqrt{n\lambda_{kn}/(p+n)}-w}\bigg]-\frac{n-p+2j}{w(p+n)}.\]
as well as
\[S_n(z)=\frac{\sqrt{n(p+n)}}{2p\sqrt{z}}\wh S_n(w_n)+\frac{n-p}{2zp},\quad w_n=\sqrt{nz/(p+n)}.\]
We now apply Proposition 12 from \cite{mp} (see also (41) in \cite{mp})
\begin{align*}
\p(|S_n(z)-&\e S_n(z)|>\varepsilon)=\p\big(\big|\wh S_n(w_n)-\e \wh S_n(w_n)\big|>2\varepsilon p\sqrt{|z|}/\sqrt{n(p+n)}\big)\\
\leqslant& \exp\Big\{-\frac{(p+n)|\Im(w_n)|^2}{2560 (m+1)}\cdot\frac{4\varepsilon^2p^2|z|}{n(p+n)}\Big\}\\
&=  \exp\Big\{-\frac{n|\Im (\sqrt{z})|^2}{2560 (m+1)}\cdot\frac{4\varepsilon^2p^2|z|}{n(p+n)}\Big\}\\&= \exp\Big\{-\frac{K n}{m+1}\Big\}
\end{align*}
with 
\[K=K(n)=\frac{4\varepsilon^2y^2|z\Im (\sqrt{z})|^2}{2560 (y+1)}+o(1),\quad n\to\infty.\]
If $m=m(n)=o(n/\log n)$ as $n\to\infty$, the latter and the Borel-Cantelli lemma imply that   $S_n(z)-\e S_n(z)\to 0$ a.s. Analogously, we obtain that 
$s_n(z)-\e s_n(z)\to 0$ a.s.

Then we proceed using the same arguments as in the proof of Theorem \ref{p1} till the definition of $I_{kn}$, i.e. 
\[I_{kn} =-\frac{ \y_{pk}^\top A_{kn}^{2}\y_{pk}/n}{1+ \y_{pk}^\top A_{kn}\y_{pk}/n}+\frac{ \z_{pk}^\top A_{kn}^2\z_{pk}/n}{1+\z_{pk}^\top A_{kn}\z_{pk}/n},\] where $A_{kn}=A_{kn}(z)=(C_{kn}-zI_p)^{-1}$, $1\leqslant k\leqslant n.$
As in the proof of Theorem \ref{p1}, we finish the proof if we show that 
\begin{equation}\label{II}
\frac{1}{p}\sum_{k=1}^n \e|I_{kn}|\to 0.
\end{equation}
  
Applying \eqref{minsum}, we arrive at
 \[\e|I_{kn}|=\e(|I_{kn}|\wedge(2/v))\leqslant\e(|\Delta_{kn}^m|\wedge(2/v))+\e(|\widehat \Delta_{kn}^m|\wedge(2/v)),\]
where
\[\Delta_{kn}^m=\frac{ \y_{pk}^\top A_{kn}^2\y_{pk}/n}{1+ \y_{pk}^\top A_{kn}\y_{pk}/n}-\frac{\tr(\Sigma_{pk} (A_{kn}^m)^2)/n}{1+\tr(\Sigma_{pk} A_{kn}^m)/n},\]
\[\widehat \Delta_{kn}^m=\frac{ \z_{pk}^\top A_{kn}^2\z_{pk}/n}{1+ \z_{pk}^\top A_{kn}\z_{pk}/n}-\frac{\tr(\Sigma_{pk} (A_{kn}^m)^2)/n}{1+\tr(\Sigma_{pk} A_{kn}^m)/n}.\] Here
$A_{kn}^m=A_{kn}^m(z)=(C_{kn}^m-zI_p)^{-1}$ and 
\[C_{kn}^m=C_{kn}-\frac{1}{n}\sum_{i=k-2m+1}^{k-1} \y_{pi}\y_{pi}^\top-\frac{1}{n}\sum_{i=k+1}^{k+2m-1} \z_{pi}\z_{pi}^\top,\quad 1\leqslant k\leqslant n,\] 
where, for simplicity, we let $\y_{pi}$ and $\z_{pi}$ be zero vectors for all $i$ such that
\[ -2m+2\leqslant i\leqslant 0\quad\text{or}\quad n+1\leqslant i\leqslant n+2m-1.\]
Note that  $C_{kn}^m$ is independent of $(\y_{pk},\z_{pk})$ for each $k=1,\ldots,n$ because of the $m$-dependence in $\{(\y_{pk},\z_{pk})\}_{k=-2m+2}^{n+2m-1}$. 

It is shown in the proof of Theorem \ref{p1} that, for any $\varepsilon>0$ and $z\in\mathbb C^+,$
\[\frac{1}{p}\sum_{k=1}^n\e(|\Delta_{kn}^m|\wedge(2/v))\leqslant C(v/|z|,1/v)\varepsilon+\frac{2}{vp}\sum_{k=1}^n\p(E_{kn}^m),\]
where $C=C(a,b)>0$ is given in Lemma \ref{l3}, $v=\Im(z)>0$ and events $E_{kn}^m=E_{kn}^m(\varepsilon)$ are defined by
\[E_{kn}^m=\bigcup_{j=1}^2\{|\y_{pk}^\top (A_{kn})^j \y_{pk}-\tr(\Sigma_{pk} (A_{kn}^m)^j)|> \varepsilon p\}.\]
A similar bound is valid when $\e(|\Delta_{kn}^m|\wedge(2/v))$ is replaced by $\e(|\wh\Delta_{kn}^m|\wedge(2/v))$, i.e. 
\[\frac{1}{p}\sum_{k=1}^n\e(|\wh\Delta_{kn}^m|\wedge(2/v))\leqslant C(v/|z|,1/v)\varepsilon+\frac{2}{vp}\sum_{k=1}^n\p\big(\wh E_{kn}^m\big),\]
where
\[\wh E_{kn}^m=\wh E_{kn}^m(\varepsilon)=\bigcup_{j=1}^2\{|\z_{pk}^\top (A_{kn})^j \z_{pk}-\tr(\Sigma_{pk} (A_{kn}^m)^j)|> \varepsilon p\}.\]

Thus, to verify \eqref{II}, we need to prove that (for any fixed $\varepsilon>0$)
\[\frac{1}{p}\sum_{k=1}^n[\p( E_{kn}^m)+\p(\wh E_{kn}^m)]\to0.\]
Fix any $\varepsilon>0.$ Obviously, $E_{kn}^m\subseteq F_{kn}^m\cup G_{kn}^m$ for 
\[F_{kn}^m=\bigcup_{j=1}^2\{|\y_{pk}^\top (A_{kn})^j \y_{pk}-\y_{pk}^\top (A_{kn}^m)^j \y_{pk}|> \varepsilon p/2\},\]
\[G_{kn}^m=\bigcup_{j=1}^2\{|\y_{pk}^\top (A_{kn}^m)^j \y_{pk}-\tr(\Sigma_{pk} (A_{kn}^m)^j)|> \varepsilon p/2\}.\]
 \begin{lemma} \label{l8} If (A3) holds, then, for any $\varepsilon,M>0,$ 
\[\frac{1}{n}\sum_{k=1}^n \sup_{A_p}\p(|\Y_{pk}^\top A_{p}\Y_{pk}-\tr(\Sigma_{pk} A_{p})|>\varepsilon p)\to 0,\quad n\to\infty,\] where all suprema are taken over all complex $p\times p$  matrices $A_p$ with $\|A_p\|\leqslant M.$
 \end{lemma} 
Using independence of $\y_{pk}$ and $(A_{kn}^m)^j$, we infer from Lemmas \ref{l0} and \ref{l8} that 
\[\frac{1}{p}\sum_{k=1}^n\p\big( G_{kn}^m\big)\leqslant\frac{2}{p}\sum_{k=1}^n \sup_{A_p}\p(|\Y_{pk}^\top A_{p}\Y_{pk}-\tr(\Sigma_{pk} A_{p})|>\varepsilon p/2)\to 0,\quad n\to\infty,\] where each suprema is taken over all complex $p\times p$  matrices  $A_p$ with spectral norms $\|A_p\|\leqslant M=1/(v\wedge v^{2})$ and we have used the fact that $p/n\to y>0.$

Similarly, $\wh E_{kn}^m\subseteq \wh F_{kn}^m\cup \wh G_{kn}^m$ for 
\[\wh F_{kn}^m=\bigcup_{j=1}^2\{|\z_{pk}^\top (A_{kn})^j \z_{pk}-\z_{pk}^\top (A_{kn}^m)^j \z_{pk}|> \varepsilon p/2\},\]
\[\wh G_{kn}^m=\bigcup_{j=1}^2\{|\z_{pk}^\top (A_{kn}^m)^j \z_{pk}-\tr(\Sigma_{pk} (A_{kn}^m)^j)|> \varepsilon p/2\}.\]
By independence of $\z_{pk}$ and $(A_{kn}^m)^j$, we derive from Lemma \ref{l0} that 
\[\frac{1}{p}\sum_{k=1}^n\p\big(\wh G_{kn}^m\big)\leqslant\frac{2}{p}\sum_{k=1}^n \sup_{A_p}\p(|\z_{pk}^\top 
 A_{p}\z_{pk}-\tr(\Sigma_{pk} A_{p})|>\varepsilon p/2),\]
where each suprema is taken over all complex $p\times p$  matrices  $A_p$ with spectral norms $\|A_p\|\leqslant M=1/(v\wedge v^{2}).$  It is shown in the proof of Proposition \ref{pz} that 
  \[\p(|\z_{pk}^\top A_{p}\z_{pk}-\tr(\Sigma_{pk} A_{p})|>\varepsilon p/2)\leqslant \e\min\Big\{\frac{16M^2\tr(\Sigma_{pk}^2)}{(\varepsilon p/2)^2},1\Big\}.\]
To show that 
\[\frac{1}{p}\sum_{k=1}^n\p\big(\wh G_{kn}^m\big)\to 0,\]  
we only need to note that (A4) implies that, for any $\varepsilon,M>0,$  
\begin{equation}\label{minS}
\frac{1}{p}\sum_{k=1}^n\e\min\Big\{\frac{16M^2\tr(\Sigma_{pk}^2)}{(\varepsilon p/2)^2},1\Big\}\to 0.
\end{equation}
Indeed, by (A4), for any $\gamma>0,$
\begin{align*}
\frac{1}{p}\sum_{k=1}^n\e\min\Big\{\frac{16M^2\tr(\Sigma_{pk}^2)}{(\varepsilon p/2)^2},1\Big\}\leqslant&\frac{n}{p}\frac{64M^2\gamma}{\varepsilon^2}+\frac{1}{p}\sum_{k=1}^n\p(\tr(\Sigma_{pk}^2)>\gamma p^2)\\
&=\frac{64M^2\gamma}{y\varepsilon^2}+o(1).\end{align*}
This clearly implies \eqref{minS}. Additionally, 
by the Markov inequality,
\[\frac{1}{p}\sum_{k=1}^n[\p(F_{kn}^m)+\p(\wh F_{kn}^m)]\leqslant \frac{2}{ \varepsilon p^2}\sum_{k=1}^n
J_k\]
with 
\[J_k=\sum_{j=1}^2\big(\e|\y_{pk}^\top (A_{kn})^j \y_{pk}-\y_{pk}^\top (A_{kn}^m)^j \y_{pk}|+\e|\z_{pk}^\top (A_{kn})^j \z_{pk}-\z_{pk}^\top (A_{kn}^m)^j \z_{pk}|\big).\]

To finish the proof of the theorem, we need to show that
\[\frac{1}{p^2}\sum_{k=1}^nJ_k\to0.\]
We need three additional lemmas.

\begin{lemma} \label{l5} Let $z\in\mathbb C^+$, $U$ be a real $p\times q$ matrix and  $C$ be a real symmetric positive semi-definite   $p\times p$ matrix. Then
\[\|(I_q+U^\top(C-zI_p)^{-1}U)^{-1}\|\leqslant \frac{|z|}{\Im(z)}.\]
\end{lemma}

\begin{lemma} \label{l6} Let $z\in\mathbb C^+$, $w\in\mathbb C^q$,  $U$ be a real $p\times q$ matrix and  $C$ be a real symmetric positive semi-definite   $p\times p$ matrix. Then
\[|w^\top A w|\leqslant\frac{\Im(z)+|z|}{|\Im(z)|^2} \|w\|^2,\]
where $A=(I_q+U^\top(C-zI_p)^{-1}U)^{-1}U^\top(C-zI_p)^{-2}U(I_q+U^\top(C-zI_p)^{-1}U)^{-1}.$
\end{lemma}

\begin{lemma} \label{l7} Let $z\in\mathbb C^+$, $y\in\mathbb C^p$,  $U$ be a real $p\times q$ matrix and  $C$ be a real symmetric positive semi-definite  $p\times p$ matrix.  Then
\begin{align*}\sum_{j=1}^2|y^\top(C+UU^\top-zI_p)^{-j}y-y^\top (C-zI_p)^{-j}y|\leqslant \frac{2(|z|+1)^2}{|\Im(z)|^2}\sum_{j=1}^2\|U^\top(C-zI_p)^{-j}y\|^2.
\end{align*}
\end{lemma}

Taking  $C=C_{kn}^m$ and 
$U$ to be $p\times (4m-2)$ matrix with columns \[n^{-1/2}\y_{pi},\;k-2m+1\leqslant i\leqslant k-1,\quad\text{and}\quad n^{-1/2}\z_{pi},\;k+1\leqslant i\leqslant k+2m-1,\] in Lemma \ref{l7}, we get $A_{kn}=(C+UU^\top-zI_p)^{-1}$ and $A_{kn}^m=(C- zI_p)^{-1}$ as well as
\begin{align*}
\sum_{j=1}^2\e|\y_{pk}^\top (A_{kn})^j \y_{pk} -\y_{pk}^\top (A_{kn}^m)^j \y_{pk}|\leqslant&\frac{2(|z|+1)^2}{v^2}\sum_{j=1}^2\e\|U^\top(C-zI_p)^{-j}\y_{pk}\|^2\\&=\frac{2(|z|+1)^2}{v^2n}\sum_{j=1}^2Q_{kj},
\end{align*}
where
\[Q_{kj}=\sum_{i=k-2m+1}^{k-1}\e|\y_{pi}^\top(A_{kn}^m)^j\y_{pk}|^2+\sum_{i=k+1}^{k+2m-1}\e|\z_{pi}^\top(A_{kn}^m)^j\y_{pk}|^2,\quad j=1,2.\]
By the same arguments,
\begin{align*}
\sum_{j=1}^2\e|\z_{pk}^\top (A_{kn})^j \z_{pk} -\z_{pk}^\top (A_{kn}^m)^j \z_{pk}|\leqslant\frac{2(|z|+1)^2}{v^2n}\sum_{j=1}^2R_{kj},\end{align*}
where
\[
R_{kj}=\sum_{i=k-2m+1}^{k-1}\e|\y_{pi}^\top(A_{kn}^m)^j\z_{pk}|^2+\sum_{i=k+1}^{k+2m-1}\e|\z_{pi}^\top(A_{kn}^m)^j\z_{pk}|^2,\quad j=1,2.\]

Let us estimate the sums $Q_{kj}$ and $R_{kj}$. There are two types of terms in $Q_{kj}$ and $R_{kj}$. Namely, terms 
\[\e|\y_{pi}^\top(A_{kn}^m)^j\z_{pk}|^2\quad \text{and}\quad \e|\y_{pi}^\top(A_{kn}^m)^j\y_{pk}|^2,\quad k-2m+1\leqslant i\leqslant k-m,\]
as well as
\[\e|\z_{pi}^\top(A_{kn}^m)^j\z_{pk}|^2\quad \text{and}\quad \e|\z_{pi}^\top(A_{kn}^m)^j\y_{pk}|^2,\quad k+m\leqslant i\leqslant k+2m-1,\]
have the form $\e|x_1^\top A x_2|^2$ for a $p\times p$ symmetric random matrix $A$ and random vectors $x_1,x_2$ in $\bR^p$ such that $(A,x_1)$ is independent of $x_2$. The rest terms have the form $\e|x_1^\top A x_2|^2$ with $(x_1,x_2)$ independent of $A$. In all cases, the spectral norm of $A$ is almost surely bounded by $M=\max\{v^{-1},v^{-2}\}$ (see Lemma \ref{l0}).

Let $A$ be any complex symmetric random $p\times p$ matrix that is independent of $(\y_{pj},\y_{pl},\z_{pj},\z_{pl})$, $1\leqslant j<l\leqslant n$ and $\|A\|\leqslant M$ a.s. Set $A^*=\overline{A^\top}=\overline{A}$ and $\e_s=\e(\cdot|\cF_s^p)$ for all $0\leqslant s\leqslant n-1$ for $\cF_s^p$ defined in the beginning of Section 3.  By the construction of $\z_{pl}=\Sigma_{pl}^{1/2}\w_{pl}$ and the law of iterated mathematical expectations,
\begin{align*}
\e|\y_{pj}^\top A& \z_{pl}|^2=\e \y_{pj}^\top A^* \z_{pl}\z_{pl}^\top A\y_{pj}=
\e \y_{pj}^\top A^*(\e_{j}\z_{pl}\z_{pl}^\top )A\y_{pj}=\\
=&\e \y_{pj}^\top A^*(\e_{j}\Sigma_{pl}) A\y_{pj}\leqslant  \e \|\e_{j}\Sigma_{pl}\| \|A\y_{pj}\|^2
\leqslant  M^2\e \|\e_{j}\Sigma_{pl}\| \|\y_{pj}\|^2.
\end{align*}
Analogously, $\e|\y_{pj}^\top A \y_{pl}|^2\leqslant M^2\e \|\e_{j}\y_{pl}\y_{pl}^\top\| \|\y_{pj}\|^2$ and
\begin{align*}
\e|\z_{pj}^\top A &\z_{pl}|^2\leqslant M^2\e \|\e_{j}\Sigma_{pl}\|\,\|\z_{pj}\|^2=\\
&= M^2\e \|\e_{j}\Sigma_{pl}\|\,\tr(\Sigma_{pj}^{1/2}\w_{pj}\w_{pj}^\top \Sigma_{pj}^{1/2})= M^2\e \|\e_{j}\Sigma_{pl}\|\,\tr(\Sigma_{pj}).
\end{align*}

Let $(A,\y_{pj})$ be independent of $(\y_{pl},\z_{pl})$, but $A$ and $\y_{pj}$ may be dependent. Then 
\begin{align}\label{aaa1}
\e|\y_{pj}^\top A \z_{pl}|^2=\e \y_{pj}^\top& A^* \z_{pl}\z_{pl}^\top A\y_{pj}=
\e \y_{pj}^\top A^*[\e(\Sigma_{pl}^{1/2}\w_{pl}\w_{pl}^\top\Sigma_{pl}^{1/2})] A\y_{pj}\leqslant\nonumber\\
\leqslant & M^2\|\e\Sigma_{pl}\| \e\|\y_{pj}\|^2= M^2\|\e\Sigma_{pl}\| \tr(\e\y_{pj}\y_{pj}^\top).
\end{align}
Similarly, $\e|\y_{pj}^\top A \y_{pl}|^2\leqslant M^2\|\e\y_{pl}\y_{pl}\| \e\|\y_{pj}\|^2.$

Let $(A,\z_{pl})$ be independent of $(\y_{pj},\z_{pj})$, but $A$ and $\z_{pl}$ may be dependent. Then 
\begin{align}\label{aaa2}
\e|\y_{pj}^\top A \z_{pl}|^2=\e| \z_{pl}^\top A &\y_{pj}|^2=\e \z_{pl}^\top  A^*\y_{pj}\y_{pj}^\top A \z_{pl}=\e \z_{pl}^\top A^*\e\y_{pj}\y_{pj}^\top A \z_{pl}\leqslant\nonumber \\
\leqslant & M^2\|\e\y_{pj}\y_{pj}^\top\| \e\|\z_{pj}\|^2= M^2\|\e\y_{pj}\y_{pj}^\top\| \tr(\e\Sigma_{pl}).
\end{align}
Similarly,  $\e|\z_{pj}^\top A \z_{pl}|^2=\e| \z_{pl}^\top A \z_{pj}|^2\leqslant M^2 \|\e\Sigma_{pj}\| \tr(\e\Sigma_{pl})$.

Combining obtained estimates, we have the following estimate 
\[Q_{k1}+Q_{k2}+R_{k1}+R_{k2}\leqslant 2M^2(S_{k1}+S_{k2}+S_{k3}+S_{k4}),\] where
\begin{align*}S_{k1}&=\sum_{i: \,m\leqslant k-i<2m} (\|\e\Sigma_{pk}\|+\|\e\y_{pk}\y_{pk}\|) \e\|\y_{pi}\|^2\\
S_{k2}&=\sum_{i:\,1\leqslant k-i<m} \e(\|\e_{i}\Sigma_{pk}\|+\|\e_{i}\y_{pk}\y_{pk}\|) \|\y_{pi}\|^2\\
S_{k3}&=\sum_{i:\,1\leqslant i-k<m} \e\|\e_{k}\Sigma_{pi}\|[\|\y_{pk}\|^2+\tr(\Sigma_{pk})]\\
S_{k4}&=\sum_{i:\,m\leqslant k-i<2m}  (\|\e\Sigma_{pk}\|+\|\e\y_{pk}\y_{pk}\|) \tr(\e \Sigma_{pi}).
\end{align*}
Define
\[T_{n1}=\sum\e(\|\e_{k}\Sigma_{pl}\|+\|\e_{k}\y_{pl}\y_{pl}^\top \|)\tr\big(\y_{pk}\y_{pk}^\top+\Sigma_{pk})\] with the sum taken over all $(k,l)$ with $1\leqslant k<l\leqslant n$ and $l<k+m$ and
\[T_{n2}=\sum(\|\e \Sigma_{pl}\|+\|\e \y_{pl}\y_{pl}^\top\|)\tr\big(\e\y_{pk}\y_{pk}^\top+\e\Sigma_{pk}\big)\]
with the sum taken over all $(k,l)$ with $1\leqslant k,l\leqslant n$ and $m\leqslant |k-l|<2m.$

The above estimates imply that
\begin{align*} \frac1{p^2}\sum_{k=1}^n J_k\leqslant \frac{4(|z|+1)^2M^2}{v^2p^2n}\sum_{k=1}^n\sum_{j=1}^2(Q_{kj}+R_{kj})\leqslant\frac{8(|z|+1)^2M^2}{v^2p^2n}(T_{n1}+T_{n2}).
\end{align*}
  Using (A5) and $n/p\to 1/y>0$, we get that
\[\frac1{p^2}\sum_{k=1}^n J_k\to 0.\]
This finishes the proof of the theorem. Q.e.d.

\noindent {\bf Proof of Proposition \ref{A5}.} Since there are $O(nm)$ terms in each sum appearing in (A5), we will verify (A5) if we show that 
\[\|\e_{l}\Sigma_{pk}\|\leqslant \e_l\|\e_{k-1}\Sigma_{pk}\|\quad\text{and}\quad \|\e_{l}\y_{pk}\y_{pk}^\top\|\leqslant \e_l\|\e_{k-1}\y_{pk}\y_{pk}^\top\|\]
for any $0\leqslant l<k\leqslant n$. All these inequalities follow from Jensen's inequality and the convexity of the spectral norm (alternatively, see E.1.Theorem in \cite{AMO}). Q.e.d.

\noindent {\bf Proof of Theorem \ref{p3}.} By the definition of $L_n$, $\wh\bY_{pn}=\bY_{pn}L_n^\top$ and $\wh\Z_{pn}=\Z_{pn}L_n^\top $ are  random matrices with $m$-dependent columns
\[\sum_{j:\, k-m<j\leqslant k}l_{kj}\y_{pj}\quad\text{and}\quad \sum_{j:\, k-m<j\leqslant k}l_{kj}\z_{pj},\quad k=1,\ldots,n,\]
where $\Y_{p1},\ldots,\Y_{pn}$ and $\bZ_{p1},\ldots,\bZ_{pn}$ are columns of $\bY_{pn}$ and  $\Z_{pn}$ correspondingly.
Therefore, we can proceed along the proof of Theorem \ref{p2} to show that  \[S_n(z)-\e S_n(z)\to 0\quad\text{and}\quad s_n(z)-\e s_n(z)\to 0\quad\text{a.s.},\quad n\to\infty,\]  where
\[S_n(z)=p^{-1}\tr\big(n^{-1}\bY_{pn}L_n^\top L_n\bY_{pn}^\top-zI_p\big)^{-1},\] \[s_n(z)=p^{-1}\tr\big(n^{-1}\Z_{pn}L_n^\top L_n\Z_{pn}^\top-zI_p\big)^{-1}.\] 

Next we will use  Lindeberg's method as in the proofs of  Theorem \ref{p1} and \ref{p2} to prove that $\e S_n(z)-\e s_n(z)\to 0$ as $n\to\infty$. Assume w.l.o.g. that  $\{(\y_{pk},\Sigma_{pk})\}_{k=1}^n$ are independent copies of $(\y_{p},\Sigma_{p})$ as well as $\z_{pk}=\Sigma_{pk}^{1/2}\w_{pk}$ for all $1\leqslant k\leqslant n$, where $\{\w_{pk}\}_{k=1}^n$ are independent standard normal vectors in $\bR^p$ that are also independent of  $\{(\y_{pk},\Sigma_{pk})\}_{k=1}^n$.

We have $|S_n(z)-s_n(z)|\leqslant \sum_{k=1}^n|I_{kn}|/p,$ where 
\begin{align*}
I_{kn} =\tr\big(n^{-1}\bY_{pn}^kL_n^\top L_n(\bY_{pn}^k)^\top-z I_p\big)^{-1}-\tr\big(n^{-1}\bY_{pn}^{k-1}L_n^\top L_n(\bY_{pn}^{k-1})^\top-z I_p\big)^{-1}
\end{align*}
for $p\times n$ matrices $\bY_{pn}^k$ defined as follows: $\bY_{pn}^0=\Z_{pn}$, $\bY_{pn}^n=\bY_{pn},$ and $\bY_{pn}^k$ is a matrix with columns 
\[\y_{pj},\quad 1\leqslant j\leqslant k,\quad\text{and}\quad\z_{pj},\quad k<j\leqslant n,\] for all $1\leqslant k<n.$

Let further $\y_{pj}$  and $\z_{pj}$ be zero vectors and $l_{kj}=0$ whenever $k,j$ do not belong to the set $\{1,\ldots,n\}$. 
Since
\[A L_n^\top L_n A^\top=\sum_{k=1}^n 
\Big(\sum_{j:\, k-m<j\leqslant k}l_{kj}a_j\Big)\Big(\sum_{j:\, k-m<j\leqslant k} l_{kj} a_j\Big)^\top.\]
for any $p\times n$ matrix $A$ with columns $a_1,\ldots,a_n$, there are symmetric positive semi-definite $p\times p$ matrices $C_{kn}$, $1\leqslant k\leqslant n$, such that 
\[\frac{1}{n}\bY_{pn}^{k-1}L_n^\top L_n(\bY_{pn}^{k-1})^\top=C_{kn}+\frac{1}{n}\sum_{s=k}^{k+m-1} \z_{ps}^{k}(\z_{ps}^{k})^\top, \]
\[\frac{1}{n}\bY_{pn}^{k}L_n^\top L_n(\bY_{pn}^{k})^\top=C_{kn}+\frac{1}{n}\sum_{s=k}^{k+m-1} \y_{ps}^{k}(\y_{ps}^{k})^\top,\]
where
\[\z_{ps}^{k}=\sum_{j:\, s-m<j<k} l_{sj}\y_{pj}+l_{sk}\z_{pk}+\sum_{j:\, k< j\leqslant s} l_{sj} \z_{pj},\]
\[\y_{ps}^{k}=\sum_{j:\, s-m<j< k} l_{sj}\y_{pj}+l_{sk}\y_{pk}+\sum_{j:\, k<j\leqslant s} l_{sj} \z_{pj} .\] By this definition, $C_{kn}$ is independent of $(\y_{pk},\z_{pk})$ for any fixed $k=1,\ldots,n.$

Denote by $U_{kn}$ and $V_{kn}$  random $p\times m$ matrices with columns  \[n^{-1/2}\z_{ps}^{k},\;k\leqslant s\leqslant k+m-1,\quad\text{and}\quad n^{-1/2}\y_{ps}^{k},\;k\leqslant s\leqslant k+m-1,\] correspondingly.  Then
\[\tr\big(n^{-1}\bY_{pn}^kL_n^\top L_n(\bY_{pn}^k)^\top-z I_p\big)^{-1}=\tr\big( C_{kn}+V_{kn} V_{kn}^\top-z I_p\big)^{-1},\]
\[\tr\big(n^{-1}\bY_{pn}^{k-1}L_n^\top L_n(\bY_{pn}^{k-1})^\top-z I_p\big)^{-1}=\tr\big( C_{kn}+U_{kn} U_{kn}^\top-z I_p\big)^{-1}.\]

By  the Sherman-Morrison-Woodbury formula,
\begin{align*} \tr&(C+UU^\top-zI_p)^{-1}=\\
&=\tr(C-zI_p)^{-1}-\tr\big((C-zI_p)^{-1}U(I_q+U^\top (C-zI_p)^{-1}U)^{-1}U^\top(C-zI_p)^{-1}\big)\\
&=\tr(C-zI_p)^{-1}-\tr\big(U^\top(C-zI_p)^{-2}U(I_q+U^\top (C-zI_p)^{-1}U)^{-1}\big)\end{align*} for any real symmetric  $p\times p$ matrix $C$ and all real $p\times q$ matrices $U$. Hence, adding and subtracting $\tr(C_{kn}-z I_p)^{-1}$ to $I_{kn}$ yield
\begin{align*} I_{kn} =&\tr\big(U_{kn}^\top A_{kn}^{2}U_{kn}(I_m+U_{kn}^\top A_{kn}U_{kn})^{-1}\big)-\tr\big(V_{kn}^\top A_{kn}^2V_{kn}(I_m+V_{kn}^\top A_{kn}V_{kn})^{-1}),\end{align*} where we set $A_{kn}=A_{kn}(z)=(C_{kn}-zI_p)^{-1}$, $1\leqslant k\leqslant n.$
\begin{lemma}\label{l10}
  Let $z\in\mathbb C^+,$ $U$ be a real $p\times q$ matrix and $C$ be a real symmetric positive semi-definite $p\times p$ matrix. Then 
\begin{align*}|\tr\big(U^\top A^{2}U(I_q+U^\top AU)^{-1}\big)|\leqslant \frac{q(|z|+\Im(z))}{|z|\Im(z)} ,\end{align*}
where  $A=(C-zI_p)^{-1}$. 
\end{lemma}

\begin{lemma} \label{l9} 
  Let $z\in\mathbb C^+,$ $U$ and $V$ be real $p\times q$ matrices and $C$ be a real symmetric positive semi-definite $p\times p$ matrix. Then 
\begin{align*}|\tr&\big(U^\top A^{2}U(I_q+U^\top AU)^{-1}\big)-\tr\big(V^\top A^{2}V(I_q+V^\top AV)^{-1}\big)|\leqslant\\&\leqslant K\sum_{j=1}^2\|U^\top A^{j}U-V^\top A^{j}V\|+K\|U^\top AU-V^\top AV\|\|U^\top AA^* U-V^\top AA^*V\|^{1/2},
\end{align*}
where $K=K(q,z)=q(|z|+1)^{3/2}/|\Im(z)|^2$, $A=(C-zI_p)^{-1}$ and $A^*=(C-\overline{z}I_p)^{-1}.$
\end{lemma}

In view of Lemma \ref{l10}, we finish the proof if we show that
\begin{equation}\label{Ip3}
\frac{1}{p}\sum_{k=1}^n\e|I_{kn}|=\frac{1}{p}\sum_{k=1}^n\e(|I_{kn}|\wedge a)\to 0,\quad n\to\infty,\end{equation}
for  $a=2m(|z|+v)/(v|z|)$ and $v=\Im(z)>0$. Fix any $\varepsilon>0$ and define 
\[D_{kn}=\bigcap_{j=1}^2\{\|U_{kn}^\top A_{kn}^{j}U_{kn}-V_{kn}^\top A_{kn}^{j}V_{kn}\|\leqslant \varepsilon\}\cap \{\|U_{kn}^\top A_{kn}A_{kn}^*U_{kn}-V_{kn}^\top A_{kn}A_{kn}^*V_{kn}\|\leqslant \varepsilon\}.\]
By Lemma \ref{l9},
\begin{align*}\frac{1}{p}\sum_{k=1}^n\e(|I_{kn}|\wedge a)=& \frac{1}{p}\sum_{k=1}^n\e(|I_{kn}|\wedge a)I(D_{kn})+\frac{1}{p}\sum_{k=1}^n\e(|I_{kn}|\wedge a)I\big(\overline{D}_{kn}\big)\\\leqslant& \frac{1}{p}\sum_{k=1}^n\e|I_{kn}|I(D_{kn})+\frac{a}{p}\sum_{k=1}^n\p\big(\overline{D}_{kn}\big)\\
\leqslant& \frac{n}pK(m,z)(2\varepsilon+\varepsilon^{3/2})+\frac{a}{p}\sum_{k=1}^n\p\big(\overline{D}_{kn}\big).
\end{align*}
Since $n/p\to1/y>0,$ to prove \eqref{Ip3} we need to show that, for any fixed $\varepsilon>0,$
\[\frac{1}{p}\sum_{k=1}^n\p\big(\overline{D}_{kn}\big)\to 0.\]
Due to Lemma \ref{l0}, it is clear that $\|A_{kn}\|\leqslant M$, $\|A_{kn}^2\|\leqslant M$
 and $\| A_{kn}A_{kn}^*\|\leqslant M$ a.s. for $M=\max\{v^{-1},v^{-2}\}.$ Moreover, $A_{kn},$ $A_{kn}^2$ and $A_{kn}A_{kn}^*=((C_{kn}-uI_p)^2+v^2I_p)^{-1}$ are symmetric matrices. As a result,
\[\p\big(\overline{D}_{kn})\leqslant 3\sup_{B_{pk}}\p(\|U_{kn}^\top B_{pk}U_{kn}-V_{kn}^\top B_{pk}V_{kn}\|>\varepsilon),\] 
where the supremum is taken over all complex symmetric $p\times p$ random matrices $B_{pk}$  such that $\|B_{pk}\|\leqslant M$ a.s. and $B_{pk}$ (with $(\y_{pj},\z_{pj})$, $j\neq k$) is independent of $(\y_{pk},\z_{pk})$.
  
  Recalling the definitions of $U_{kn}$ and $V_{kn}$, we  write $U_{kn}=\wh U_{kn} L_{kn}^\top$ and
$V_{kn}=\wh V_{kn} L_{kn}^\top$, where
\[L_{kn}=\begin{pmatrix}
l_{k,k-m+1}&\ldots&l_{kk}&0&\ldots&0\\
0&l_{k+1,k-m+1}&\ldots &l_{k+1,k+1}&0&\ldots\\
&&&\ldots&&
\\
\ldots&0&l_{k+m-1,k}&\ldots& &l_{k+m-1,k+m-1}\\
\end{pmatrix}\] is real $m\times (2m-1)$  matrix, $ \wh U_{kn}$ is a $p\times (2m-1)$ matrix with columns  
\[n^{-1/2}\y_{pj},\;k-m+1\leqslant j< k,\quad\text{and} \quad n^{-1/2}\z_{pj},\;k\leqslant j\leqslant k+m-1,\] and  $ \wh V_{kn}$ is a $p\times (2m-1)$ matrix with columns  
\[n^{-1/2}\y_{pj},\;k-m+1\leqslant j\leqslant k,\quad\text{and} \quad n^{-1/2}\z_{pj},\;k< j\leqslant k+m-1.\]

By the assumption of Theorem \ref{p3}, entries of $L_n$ are uniformly bounded over $n$. In addition, $L_{kn}$, $1\leqslant k\leqslant n$, is a submatrix of $L_n$ which size does not depend on $n$ and is equal to $m\times(2m-1)$. Therefore, 
\[ S=\sup\{\|L_{kn}\|\|L_{kn}^\top\|:n\geqslant 1,k=1,\ldots,n\}<\infty.\]
 Hence,
\begin{align*}
\|U_{kn}^\top B_{pk}U_{kn}-V_{kn}^\top B_{pk}V_{kn}\|=&
\|L_{kn}(\wh U_{kn}^\top B_{pk}\wh U_{kn}-\wh V_{kn}^\top B_{pk}\wh V_{kn})L_{kn}^\top\|\\
\leqslant& S\|\wh U_{kn}^\top B_{pk}\wh U_{kn}-\wh V_{kn}^\top B_{pk}\wh V_{kn}\|.
\end{align*}
For any $p\times p$ matrix $A$ and any $p\times (2m-1)$ matrix $U$ with columns $u_1,\ldots,u_{2m-1}$, 
\[U^\top AU=(u_i^\top Au_j)_{i,j=1}^{2m-1}.\]
Hence, applying the bound 
\[\|Q\|^2\leqslant \|Q^*Q\|\leqslant\tr(Q^*Q)=\sum_{i,j=1}^{d}|q_{ij}|^2\]
valid for any complex $d\times d$ matrix $Q=(q_{ij})_{i,j=1}^{d},$ we get
\begin{align*}
\|\wh U_{kn}^\top B_{pk}\wh U_{kn}-\wh V_{kn}^\top B_{pk}\wh V_{kn}\|&\leqslant \big(J_{kn}+|\y_{pk}^\top B_{pk}\y_{pk}-\z_{pk}^\top B_{pk}\z_{pk}|^2/n^2\big)^{1/2}\\
 &\leqslant \sqrt{J_{kn}}+|\y_{pk}^\top B_{pk}\y_{pk}-\z_{pk}^\top B_{pk}\z_{pk}|/n
\end{align*}
where
\[J_{kn}=\frac{2}{n^2}\sum_{i=k-m+1}^{k-1}|\y_{pi}^\top B_{pk}\y_{pk}-\y_{pi}^\top B_{pk}\z_{pk}|^2+\frac{2}{n^2}\sum_{i=k+1}^{k+m-1}|\z_{pi}^\top B_{pk}\y_{pk}-\z_{pi}^\top B_{pk}\z_{pk}|^2.\]
Gathering together these bounds,
\begin{align*}
\p(\|U_{kn}^\top B_{pk}U_{kn}-&V_{kn}^\top B_{pk}V_{kn}\|>\varepsilon)\leqslant 
\p(\|\wh U_{kn}^\top B_{pk}\wh U_{kn}-\wh V_{kn}^\top B_{pk}\wh V_{kn}\|>\varepsilon/S)\\
\leqslant&
\p(J_{kn}^{1/2}>\varepsilon/(3S)\big)+\p(|\z_{pk}^\top B_{pk}\z_{pk}-\tr(\Sigma_{pk}B_{pk})|>\varepsilon/(3S))\\
&+\p(|\y_{pk}^\top B_{pk}\y_{pk}-\tr(\Sigma_{pk}B_{pk})|>\varepsilon/(3S))\\
\leqslant &\frac{9S^2}{\varepsilon^2}\,\e J_{kn}+\sup_{A_p}\p(|\z_{p}^\top A_{p}\z_{p}-\tr(\Sigma_{p}A_{p})|>\varepsilon/(3S))\\
&+\sup_{A_p}\p(|\y_{p}^\top A_{p}\y_{p}-\tr(\Sigma_{p}A_{p})|>\varepsilon/(3S)),
\end{align*}
where the supremum is taken over all complex $p\times p$ matrices $A_p$ with $\|A_p\|\leqslant M.$

By  (A1) and Lemma \ref{l2}, 
\[\sup_{A_p}\p(|\y_{p}^\top A_{p}\y_{p}-\tr(\Sigma_{p}A_{p})|>\varepsilon/(3S))=o(1),\quad p\to\infty.\]
By (A2) and Proposition \ref{pz},
\[\sup_{A_p}\p(|\z_{p}^\top A_{p}\z_{p}-\tr(\Sigma_{p}A_{p})|>\varepsilon/(3S))=o(1).\]
All these bounds guarantee that (for fixed $\varepsilon>0$ and $D_{kn}=D_{kn}(\varepsilon)$)
\[
\frac{1}{p}\sum_{k=1}^n\p\big(\overline{D}_{kn}\big)\to 0
\]
whenever $m$ is fixed, $n\to\infty$, $p/n\to y>0$, and  
\[ \frac{1}{p}\sum_{k=1}^n\sup_{B_{pk}}\e J_{kn}\to 0.\]
Thus we need to verify the last relation in order to finish the proof of the theorem.

By the Cauchy inequality,
\begin{align*}
\e J_{kn}\leqslant&
\frac{4}{n^2}\sum_{i=k-m+1}^{k-1}\big(\e|\y_{pi}^\top B_{pk}\y_{pk}|^2+\e|\y_{pi}^\top B_{pk}\z_{pk}|^2\big)\\&+\frac{4}{n^2}\sum_{i=k+1}^{k+m-1}\big(\e|\z_{pi}^\top B_{pk}\y_{pk}|^2+\e|\z_{pi}^\top B_{pk}\z_{pk}|^2\big).
\end{align*}
Using independence of $(\y_{pk},\z_{pk})$ and $(\y_{pi},\z_{pi}, B_{pk})$, $i\neq k,$ inequality $\|B_{pk}\|\leqslant M$ a.s. and arguing as in the end of the proof of Theorem \ref{p2} (see \eqref{aaa1} and \eqref{aaa2}),
\begin{align*}
\e|\y_{pi}^\top B_{pk}\y_{pk}|^2&\leqslant M^2\|\e\y_{pk}\y_{pk}^\top\|\tr(\e\y_{pi}\y_{pi}^\top)=
M^2\|\e\y_{p}\y_{p}^\top\|\tr(\e\y_{p}\y_{p}^\top),\\
\e|\z_{pi}^\top B_{pk}\y_{pk}|^2&\leqslant M^2\|\e\y_{pk}\y_{pk}^\top\|\tr(\e\Sigma_{pi})=
M^2\|\e\y_{p}\y_{p}^\top\|\tr(\e\Sigma_{p}),\\
\e|\y_{pi}^\top B_{pk}\z_{pk}|^2&\leqslant M^2\|\e\Sigma_{pk}\|\tr(\e\y_{pi}\y_{pi}^\top)=
M^2\|\e\Sigma_{p}\|\tr(\e\y_{p}\y_{p}^\top),\\
\e|\z_{pi}^\top B_{pk}\z_{pk}|^2&\leqslant M^2\|\e\Sigma_{pk}\|\tr(\e\Sigma_{pi})=M^2\|\e\Sigma_{p}\|\tr(\e\Sigma_{p}).
\end{align*}
It now follows from (A6) and \eqref{Dnorm} that 
\[\frac{1}{p}\sum_{k=1}^n\sup_{B_{pk}}\e J_{kn}\leqslant \frac{4M^2(m-1)}{pn}\big(\|\e\Sigma_{p}\|+\|\e\y_{p}\y_{p}^\top\|\big)\tr\big(\e\Sigma_{p}+\e\y_{p}\y_{p}^\top\big)=o(1).\]
This finishes the proof of the theorem. Q.e.d.

\noindent {\bf Proof of Theorem \ref{t1}.}  We will prove the first inequality by the same arguments as in the proof of Theorem 3 in \cite{G72}. Write $a=(a_1,\ldots,a_p).$  By Lemma 1 in \cite{G72},\[|(\X_p^\top a)^4-24 T|\leqslant C_0\sum_{j=0}^2|\X_p^\top a|^j|S|^{4-j},\]
where $C_0>0$ is a universal constant,
\[T=\sum_{i<j<k<l}a_ia_ja_ka_lX_iX_jX_kX_l,\quad S=\Big(\sum_{i=1}^pa_i^2X_i^2\Big)^{1/2},\]
hereinafter $i,j,k,l$ are any numbers in $\{1,\ldots,p\}$.
Hence, by Minkowski's inequality,
\begin{align*}
\e(\X_p^\top a)^4\leqslant &24 \e T+C_0\sum_{j=0}^2\e |\X_p^\top a|^j|S|^{4-j}\\
\leqslant& 
24 \e T+C_0\sum_{j=0}^2[\e(\X_p^\top a)^4]^{j/4}(\e S^4)^{1-j/4}.
\end{align*}
 By  the Cauchy-Schwartz and \eqref{vphi},
\[
|\e T|\leqslant\sum_{ i<j<k<l}|a_ia_ja_ka_l| \min\{\varphi_{j-i},\varphi_{k-j},\varphi_{l-k}\}\leqslant \sqrt{J_1J_2}\]
with 
\[J_1=\sum_{ i<j<k<l}a_i^2a_j^2 \min\{\varphi_{k-j},\varphi_{l-k}\},\quad J_2=\sum_{ i<j<k<l}a_k^2a_l^2 \min\{\varphi_{j-i},\varphi_{k-j}\}.\]
In addition, \begin{align*}
J_1\leqslant&\sum_{i<j}a_{i}^2a_j^2\sum_{k:\,k>j}\Big((k-j)\varphi_{k-j}+\sum_{l:\,l-k>k-j}\varphi_{l-k}\Big) \leqslant\\
& \quad\leqslant \frac{\|a\|^4}2\,\Phi_1+\frac{\|a\|^4}2\sum_{q=1}^\infty \sum_{r=q+1}^\infty\varphi_r=\\
&\qquad=\frac{\|a\|^4}2\,\Phi_1+\frac{\|a\|^4}2\sum_{r=2}^\infty (r-1)\varphi_r\leqslant \Phi_1 \|a\|^4.
\end{align*}
Similar arguments yield 
\begin{align*}
J_2\leqslant&
\sum_{k<l}
a_{k}^2a_l^2\sum_{j:\,j<k}\Big((k-j)\varphi_{k-j}+\sum_{i:\,j-i>k-j}\varphi_{j-i}\Big)\leqslant
\Phi_1 \|a\|^4.
\end{align*}
Let us also note that 
\[\e|S|^4=\sum_{i,j=1}^pa_i^2 a_j^2\e X_i^2X_j^2\leqslant \Phi_0\|a\|^4.\]
Combining the above estimates, we infer that
\[\e(\X_p^\top a)^4\leqslant 24 (\Phi_0+\Phi_1)\|a\|^4 +C_0\sum_{j=0}^2[\e(\X_p^\top a)^4]^{j/4}( (\Phi_0+\Phi_1)\|a\|^4)^{1-j/4}.\]
Put $R=(\e(\X_p^\top a)^4)^{1/4}(\Phi_0+\Phi_1)^{-1/4}/\|a\|$. Then \[R^4\leqslant 24  +C_0+C_0R+C_0R^2.\] Hence, $R\leqslant R_0$, where $R_0>0$ is the largest root of the equation \[x^4=24  +C_0+C_0x+C_0x^2.\] Finally we conclude that
 $\e|\X^\top a|^4\leqslant R_0^4(\Phi_0+\Phi_1)\|a\|^4.$

We now verify the second inequality. Let $A=(a_{ij})_{i,j=1}^p$ and $a_{ii}=0$, $1\leqslant i\leqslant p$. Since $\X_p^\top A\X_p=\X_p^\top B\X_p$ and
\begin{equation}\label{BB}
\tr(BB^\top)=\sum_{i,j=1}^p \Big(\frac{a_{ij}+a_{ji}}2\Big)^2\leqslant
\sum_{i,j=1}^p \frac{a_{ij}^2+a_{ji}^2}2=\sum_{i,j=1}^p a_{ij}^2=\tr(AA^\top)
\end{equation}
for $B=(A^\top+A)/2$, we may assume that $A=A^\top.$  Then 
\begin{align*}
\e|\X_p^\top A\X_p|^2=4
&\e\Big|\sum_{i=1}^{p-1} X_i \sum_{k=i+1}^p a_{ik} X_k\Big|^2\\
=&4\sum_{i=1}^{p-1} \e X_i^2\Big|\sum_{k=i+1}^p a_{ik} X_k\Big|^2+8\sum_{i<j}\e X_iX_j\Big(\sum_{k=i+1}^p a_{ik} X_k\Big)\Big(\sum_{k=j+1}^p a_{jk} X_k\Big)\\
= &4I_1+8I_2+8I_3+8I_4,
\end{align*}
where
\begin{align*}
I_1&=\sum_{i=1}^{p-1} \e X_i^2\Big|\sum_{k=i+1}^p a_{ik} X_k\Big|^2,\\
I_2&=\sum_{i<j}\e X_i\Big(\sum_{k=i+1}^{j-1} a_{ik} X_k\Big) X_j\Big(\sum_{k=j+1}^p a_{jk} X_k\Big),\\
I_3&=\sum_{i<j}a_{ij}\e X_i X_j^2\Big(\sum_{k=j+1}^p a_{jk} X_k\Big),\\
I_4&=\sum_{i<j}\e X_iX_j\Big(\sum_{k=j+1}^p a_{ik} X_k\Big)\Big(\sum_{k=j+1}^p a_{jk} X_k\Big)
\end{align*}
and sums over the empty set are zeros.
\\
{\it Control of $I_1$.} By the Cauchy-Schwartz inequality and the first part of  Theorem \ref{t1},
\begin{align*}
I_1\leqslant &
\sum_{i=1}^{p-1} \sqrt{\e X_i^4}\Big(\e\Big|\sum_{k=i+1}^p a_{ik} X_k\Big|^4\Big)^{1/2}\\
\leqslant&
 C(\Phi_0+\Phi_1) \sum_{i=1}^{p-1} \sum_{k=i+1}^p a_{ik}^2\\
 &=C(\Phi_0+\Phi_1)\frac{\tr(A^2)}2
\end{align*}
\\
{\it Control of $I_2$.} By the Cauchy-Schwartz inequality and \eqref{vphi},
\begin{align*} I_2\leqslant&\sum_{i<k<j<l}|a_{ik}a_{jl}|\,|\e X_iX_kX_jX_l|\\
\leqslant&\sum_{i<k<j<l}|a_{ik}a_{jl}|\min\{\varphi_{k-i},\varphi_{j-k},\varphi_{l-j}\}\\
\leqslant& \sqrt{I_5I_6},\end{align*}
where
\[I_5=\sum_{i<k<j<l}
a_{ik}^2\min\{\varphi_{j-k},\varphi_{l-j}\},\quad 
I_6=\sum_{i<k<j<l}
a_{jl}^2\min\{\varphi_{k-i},\varphi_{j-k}\}.
\]
Additionally,
\begin{align*}
I_5\leqslant&\sum_{i<k}a_{ik}^2\sum_{j:\,j>k}\Big((j-k)\varphi_{j-k}+\sum_{l:\,l-j>j-k}\varphi_{l-j}\Big)\\
\leqslant& \frac{\tr(A^2)}2\,\Phi_1+\frac{\tr(A^2)}2\sum_{q=1}^\infty \sum_{r=q+1}^\infty\varphi_r=\\&=\frac{\tr(A^2)}2\,\Phi_1+\frac{\tr(A^2)}2\sum_{r=2}^\infty (r-1)\varphi_r\leqslant \tr(A^2)\Phi_1.
\end{align*}
We similarly derive that
\begin{align*}
I_6\leqslant&
\sum_{j<l}
a_{jl}^2\sum_{k:\,k<j}\Big((j-k)\varphi_{j-k}+\sum_{i:\,k-i>j-k}\varphi_{k-i}\Big)\leqslant
\tr(A^2)\Phi_1.
\end{align*}
Hence, $I_2\leqslant\tr(A^2)\Phi_1 $.
\\
{\it Control of $I_3$.} By the Cauchy-Schwartz inequality and the first part of  Theorem \ref{t1},
\begin{align*}
I_3=&\sum_{j=2}^{p-1}\e  X_j^2\Big(\sum_{i=1}^{j-1} a_{ij} X_i\Big)\Big(\sum_{k=j+1}^p a_{jk} X_k\Big)\\
&\leqslant\sum_{j=2}^{p-1}\sqrt{\e  X_j^4}\bigg[\e\Big(\sum_{i=1}^{j-1} a_{ij} X_i\Big)^4\e\Big(\sum_{k=j+1}^p a_{jk} X_k\Big)^4\bigg]^{1/4}\\&\leqslant\sqrt{C(\Phi_0+\Phi_1)I_7I_8},
\end{align*}
where
\[I_7=\sum_{j=2}^{p-1}\bigg[\e\Big(\sum_{i=1}^{j-1} a_{ij} X_i\Big)^4\bigg]^{1/2},\quad I_8=\sum_{j=2}^{p-1}\bigg[\e\Big(\sum_{k=j+1}^p a_{jk} X_k\Big)^4\bigg]^{1/2}.\]
By the first inequality in Theorem \ref{t1},
\[I_7\leqslant K\sum_{j=2}^{p-1}\sum_{i=1}^{j-1} a_{ij}^2\leqslant\frac{K\tr(A^2)}{2},\quad 
I_8\leqslant K\sum_{j=2}^{p-1}\sum_{k=j+1}^p a_{jk}^2\leqslant\frac{K\tr(A^2)}{2},\]
where $K=\sqrt{C(\Phi_0+\Phi_1)}.$ As a result, $I_3\leqslant C(\Phi_0+\Phi_1)\tr(A^2)/2$.
\\
{\it Control of $I_4$.} We have $I_4=I_9+I_{10}+I_{11}$, where
\[I_9=\sum_{i<j<k}\e (a_{ik }X_i) (a_{jk} X_j)X_k^2,\;\; I_{10}=
\sum_{i<j<k<l}a_{ik}a_{jl}\e X_iX_jX_{k}X_l,\]\[
I_{11}=\sum_{i<j<k<l}a_{il}a_{jk}\e X_iX_jX_{k}X_l.\]
By the first inequality in Theorem \ref{t1},
\begin{align*}
I_9=&\frac{1}{2}\sum_{k=3}^p \e\Big(\sum_{i=1}^{k-1}a_{ik}X_i\Big)^2X_k^2-
\frac{1}{2}\sum_{k=3}^p \e X_k^2\sum_{i=1}^{k-1}a_{ik}^2X_i^2\\
&\leqslant\frac{1}{2} \sum_{k=3}^p\bigg[ \e\Big(\sum_{i=1}^{k-1}a_{ik}X_i\Big)^4\bigg]^{1/2}\sqrt{\e X_k^4}\\&\leqslant C(\Phi_0+\Phi_1)
\sum_{k=3}^p \sum_{i=1}^{k-1}\frac{a_{ik}^2}2\\&\leqslant C(\Phi_0+\Phi_1)\,\frac{\tr(A^2)}{4}.
\end{align*}
We will estimate  $I_{10}$ and $I_{11}$ in the same way as $I_2.$
It follows from the Cauchy-Schwartz inequality that $I_{10}\leqslant \sqrt{I_{12}I_{13}}$ and  $I_{11}\leqslant \sqrt{I_{14}I_{15}}$ with
\[I_{12}=\sum_{i<j<k<l}a_{ik}^2 \min\{\varphi_{j-i},\varphi_{l-k}\},\quad I_{13}=\sum_{i<j<k<l}a_{jl}^2 \min\{\varphi_{j-i},\varphi_{l-k}\},\]
\[I_{14}=\sum_{i<j<k<l}a_{il}^2 \min\{\varphi_{j-i},\varphi_{l-k}\},\quad I_{15}=\sum_{i<j<k<l}a_{jk}^2 \min\{\varphi_{j-i},\varphi_{l-k}\},\]
As previously, we have 
\begin{align*}
I_{12}&\leqslant\sum_{i<k}a_{ik}^2\sum_{j:\,i<j<k}\Big((j-i)\varphi_{j-i}+\sum_{l:\,l-k>j-i}\varphi_{l-k}\Big)\leqslant \tr(A^2)\Phi_1,\\
I_{13}&\leqslant\sum_{j<l}a_{jl}^2\sum_{k:\,j<k<l}\Big((l-k)\varphi_{l-k}+\sum_{i:\,j-i>l-k}\varphi_{j-i}\Big)\leqslant \tr(A^2)\Phi_1,\\
I_{14}&\leqslant\sum_{i<l}a_{il}^2\sum_{k:\,i<k<l}\Big((l-k)\varphi_{l-k}+\sum_{j:\,j-i>l-k,\,j<k}\varphi_{j-i}\Big)\leqslant \tr(A^2)\Phi_1,\\
I_{15}&\leqslant\sum_{j<k}a_{jk}^2\sum_{i:\,i<j}\Big((j-i)\varphi_{j-i}+\sum_{l:\,l-k>j-i}\varphi_{l-k}\Big)\leqslant \tr(A^2)\Phi_1.
\end{align*}
Thus, $I_4\leqslant C(\Phi_0+\Phi_1)\tr(A^2)/4+ 2\tr(A^2)\Phi_1$.

Combining all above estimates, we get
\[\e|\X_p^\top A\X_p|^2\leqslant (8 C(\Phi_0+\Phi_1)+24\Phi_1)\tr(A^2).\]
Q.e.d.
\\
\noindent{\bf Proof of Theorem \ref{t2}.}
Let $A=(a_{ij})_{i,j=1}^p$ and $D$ be the $p\times p$ diagonal matrix with diagonal entries $a_{11},\ldots,a_{pp}$.
By Theorem \ref{t1},
 \[\e|\X_p^\top (A-D)\X_p |^2 \leqslant C(\Phi_0+\Phi_1)\tr((A-D)(A-D)^\top).\] 
In addition,  $\tr(\Sigma_p A)=\tr(\Sigma_p D)$ for the diagonal matrix $\Sigma_p$ with diagonal entries $\e X_1^2,\ldots,\e X_p^2$ as well as
 \[\e|\X_p^\top A\X_p -\tr(\Sigma_p A)|^2\leqslant 2\e|\X_p^\top D\X_p -\tr(\Sigma_p D)|^2+2\e|\X_p^\top (A-D)\X_p |^2 .\] 
Since \[\tr (AA^\top)=\tr((A-D)(A-D)^\top)+\tr(D^2),\] we only need to bound $\e|\X_p^\top D\X_p -\tr(\Sigma_p D)|^2$ from above by $\tr(D^2)$ up to a constant factor. Write $D=D_1-D_2$, where $D_i$ are diagonal matrices with non-negative diagonal elements.  Since
 \[\e|\X_p^\top D\X_p -\tr(\Sigma_p  D)|^2\leqslant 2\sum_{i=1}^2\e|\X_p^\top D_i\X_p -\tr(\Sigma_p  D_i)|^2,\] 
 we may assume w.l.o.g. that diagonal elements of $D$ are non-negative.

We see that 
\begin{align*}
\e|\X_p^\top D\X_p -\tr(\Sigma_p D)|^2=&\e\Big|\sum_{i=1}^pa_{ii}(X_i^2-\e X_i^2)\Big|^2\\
=& \sum_{i=1}^p a_{ii}^2\var( X_i^2)+\sum_{i\neq j}a_{ii}a_{jj}\cov( X_i^2, X_j^2) 
\\\leqslant&\Phi_0\tr(D^2)+\sum_{i\neq j}\frac{a_{ii}^2+a_{jj}^2}2 \phi_{|i-j|}\\
\leqslant&\Phi_0\tr(D^2)+\sum_{i=1}^p a_{ii}^2\sum_{j: j\neq i}\phi_{|i-j|}\\
\leqslant& 2\tr(D^2)\Big(\Phi_0+\sum_{k=1}^\infty \phi_k\Big)= 2(\Phi_0+\Phi_2)\tr (D^2).
\end{align*}
Combining the above bounds, we get the desired inequality. Q.e.d.\\
\\
{{\bf Proof of Theorem \ref{t5}.}}  As in the proof of Theorem \ref{t1}, we may assume that $A$ is symmetric (see \eqref{BB}).
Note that \[V_k=\Big(\sum_{i=1}^{k-1}a_{ik} X_i\Big)X_k,\quad 2\leqslant k\leqslant p,\]
is a martingale difference sequence and
\[\X_p^\top A\X_p=2\sum_{1\leqslant i<k\leqslant p}a_{ik} X_iX_k=2\sum_{k=2}^p V_k,\] we infer
\begin{align*}\e\Big|\sum_{k=2}^pV_k\Big|^2=&\sum_{k=2}^p\e V_k^2= \sum_{k=2}^p\e \Big(\sum_{i=1}^{k-1}a_{ik} X_i\Big)^2\e(X_k^2|X_1,\ldots,X_{k-1})\leqslant\\
&\leqslant M \sum_{k=2}^p\e \Big(\sum_{i=1}^{k-1}a_{ik} X_i\Big)^2=M^2\sum_{k=2}^p\sum_{i=1}^{k-1}a_{ik}^2\leqslant\frac{M^2\tr(A^2)}{2}.
\end{align*} 
The latter implies the desired bound. Q.e.d.\\
{{\bf Proof of Theorem \ref{t6}.}} Write $A=(a_{ik})_{i,k=1}^p.$  As in the proof of Theorem \ref{t1}, we may assume that $A$ is symmetric (see \eqref{BB}).   We have
\[\X_p^\top A\X_p=\sum_{k=1}^pa_{kk}X_{k}^2+\sum_{ i\neq k}a_{ik} X_iX_k\]
and, by the Cauchy-Schwartz inequality, 
\[\e|\X_p^\top A\X_p -\tr(A)|\leqslant \e\Big|\sum_{k=1}^pa_{kk}(X_{k}^2-1)\Big|+\Big(\e\Big|\sum_{ i\neq k}a_{ik} X_iX_k\Big|^2\Big)^{1/2}.\]
By the Burkholder-Davis-Gundy inequality,
\begin{align*}\e\Big|\sum_{k=1}^pa_{kk}(X_{k}^2-1)\Big|\leqslant& C \e\Big|\sum_{k=1}^pa_{kk}^2(X_{k}^2-1)^2\Big|^{1/2},
\end{align*}
where $C$ is a universal constant. Using inequality
\begin{equation}\label{xy}
\sqrt{x+y}\leqslant\sqrt{x}+\sqrt{y},\quad x,y\geqslant0,
\end{equation}
we derive that
\[\e\Big|\sum_{k=1}^pa_{kk}^2(X_{k}^2-1)^2\Big|^{1/2}\leqslant I_1+I_2,\]
where
\begin{align*}
I_1=&\e\Big|\sum_{k=1}^pa_{kk}^2(X_{k}^2-1)^2I(|X_{k}^2-1|\leqslant b^2 )\Big|^{1/2},\\
I_2=& \e\Big|\sum_{k=1}^p a_{kk}^2(X_{k}^2-1)^2I(|X_{k}^2-1|>b^2)\Big|^{1/2}.
\end{align*}
By Jensen's inequality, 
\[I_1\leqslant \Big|\sum_{k=1}^pa_{kk}^2\e(X_{k}^2-1)^2I(|X_{k}^2-1|\leqslant b^2)\Big|^{1/2}\leqslant b\sqrt{2\,\tr(AA^\top)}.\]
Here we also used that 
\[\e(X_{k}^2-1)^2I(|X_{k}^2-1|\leqslant b^2)\leqslant b^2\e|X_{k}^2-1|\leqslant 2b^2,\quad 1\leqslant k\leqslant n.\]
In addition, it follows from \eqref{xy} that 
\[I_2\leqslant\sum_{k=1}^p |a_{kk}|\e |X_{k}^2-1|I(|X_{k}^2-1|> b^2).\]
By Theorem \ref{t5},
\[\e\Big|\sum_{ i\neq k}a_{ik} X_iX_k\Big|^2\leqslant 2\tr(AA^\top ).\]
Combining the above estimates, we get the desired result. 
Q.e.d.
\\
{\bf Proof of Corollary \ref{t4}.} If $\y_p\in\mathcal X_p$, then $\Gamma_{pn} \X_n\to \y_p $ in probability and in mean square as $n\to\infty$ for some $p\times n$ matrices $\Gamma_{pn}$ and $\X_n=(X_1,\ldots,X_n)$.

Since $(X_k)_{k\geqslant 1}$ is an orthonormal sequence, we have  \[\Gamma_{pn}\Gamma_{pn}^\top=\e (\Gamma_{pn} \X_n)(\Gamma_{pn} \X_n)^\top\to \e \Y_p\Y_p^\top=\Sigma_p,	\] \[\X_n^\top(\Gamma_{pn}^\top A_p \Gamma_{pn})\X_n=(\Gamma_{pn} \X_n)^\top A_p \Gamma_{pn} \X_n\pto \y_p^\top A_p\y_p\]
and $\tr(\Gamma_{pn}^\top A_p \Gamma_{pn})=\tr( \Gamma_{pn}\Gamma_{pn}^\top A_p)\to \tr(\Sigma_pA_p)$  as $n\to\infty$.   We need the following version of Fatou's lemma:
\begin{equation}\label{FL}
\text{If $\xi_n\pto \xi$, then $\e |\xi|\leqslant \varliminf\limits_{n\to\infty}\e |\xi_n|.$}
\end{equation}
By this lemma and  Theorem \ref{t2},
\begin{align*}
\e|\Y_p^\top A_p\Y_p-\tr(\Sigma_p A_p)|^2\leqslant & \varliminf\limits_{n\to\infty} \e|\X_n^\top( \Gamma_{pn}^\top A_p\Gamma_{pn}) \X_n-\tr(\Gamma_{pn}^\top A_p\Gamma_{pn})|^2\\
\leqslant&\varliminf\limits_{n\to\infty} C(\Phi_0+\Phi_1+\Phi_2) \tr(\Gamma_{pn}^\top A_p\Gamma_{pn}\Gamma_{pn}^\top A_p\Gamma_{pn})
\end{align*}
Note that 
\[
\tr(\Gamma_{pn}^\top A_p\Gamma_{pn}\Gamma_{pn}^\top A_p\Gamma_{pn})=\tr(\Gamma_{pn}\Gamma_{pn}^\top A_p\Gamma_{pn}\Gamma_{pn}^\top A_p)\to\tr(\Sigma_pA_p\Sigma_pA_p).\]
Let $A_p^{1/2}$ be the principal square root of $A_p$. Then \[R:=\tr(\Sigma_pA_p\Sigma_pA_p)=\tr(A_p^{1/2}\Sigma_pA_p\Sigma_pA_p^{1/2}).\]
Since $\|A_p\|I_p-A_p$ is positive semi-definite, then $Q^\top(\|A_p\|I_p-A_p) Q$ is positive semi-definite for any matrix $Q$. Taking $Q=\Sigma_pA_p^{1/2}$, we get
\begin{equation}\label{QQ}
R=\tr(Q^\top A_pQ)\leqslant\|A_p\|\tr(Q^\top Q)=\|A_p\|\tr(QQ^\top)=\|A_p\|\tr(\Sigma_pA_p\Sigma_p).
\end{equation} 
Analogously, $\tr(\Sigma_pA_p\Sigma_p)\leqslant \|A_p\|\tr(\Sigma_p^2).$
 Hence, we obtain the desired bound. 
 
 In particular, we get that
 \[\e|[\y_p^\top A_p\y_p -\tr(\Sigma_pA_p)]/p|^2=O(\|A_p\|^2\tr( \Sigma_p^2)/p^2),\quad p\to\infty.\]  
Thus $(A1)$ holds when $(A2)$ holds. Q.e.d.\\
{\bf Proof of Corollary \ref{t7}.} Arguing as in the proof of Corollary \ref{t4} and using the same notation, we derive 
\begin{align*}
\e|\Y_p^\top A_p\Y_p-\tr(\Sigma_p A_p)|\leqslant & \varliminf\limits_{n\to\infty} \e|\X_n^\top( \Gamma_{pn}^\top A_p\Gamma_{pn}) \X_n-\tr(\Gamma_{pn}^\top A_p\Gamma_{pn})|
\end{align*}
If $A_p$ is a symmetric positive semi-definite matrix, then  $B_n=\Gamma_{pn}^\top A_p\Gamma_{pn}$ has the same properties. Particularly, diagonal entries of   $B_n$ are non-negative. Theorem \ref{t6} yields
\[ \e|\X_n^\top B_n \X_n-\tr(B_n)|\leqslant Cb\sqrt{\tr(B_nB_n^\top)}+CL(b)\tr(B_n)\]
for all $b,n>1$, some universal constant $C>0$ and 
\[L(b)=\sup_{k\geqslant 1}\e|X_k^2-1|I(|X_k^2-1|> b^2).\] 
As in the proof of Corollary \ref{t4}, we get
\[\lim_{n\to\infty}\tr(B_nB_n^\top)=\tr(\Sigma_pA_p\Sigma_pA_p)\leqslant\|A_p\|^2\tr(\Sigma_p^2),\]
\[\lim_{n\to\infty} \tr(B_n)=\tr(\Sigma_pA_p)\leqslant \|A_p\|\tr(\Sigma_p).\]
Thus we have proven that
\[\e|\Y_p^\top A_p\Y_p-\tr(\Sigma_p A_p)|\leqslant Cb\|A_p\|\sqrt{\tr(\Sigma_p^2)}+CL(b)\|A_p\|\tr(\Sigma_p).\]

Take any sequence of symmetric positive semi-definite $p\times p$ matrices  $A_p$ with $\|A_p\|=O(1)$  ($p=1,2,\ldots$) and assume that $\tr(\Sigma_p^2)=o(p^2)$ holds and $\tr(\Sigma_p)=O(p)$. Then, for any $b>1,$
\[\e|[\Y_p^\top A_p\Y_p-\tr(\Sigma_p A_p)]/p|\leqslant o(1)+L(b)T_p,\quad p\to\infty,\]
where $T_p=O(1)$ does not depend on $b.$  If $\{X_k^2\}_{k=1}^\infty$ is a uniformly integrable family, then $L(b)\to 0$ as $b\to\infty,$
 since \[\e|X_k^2-1|I(|X_k^2-1|> b^2)\leqslant 2\e X_k^2I(X_k^2> b^2+1)\leqslant  2\e X_k^2I(|X_k|> b)\to 0\]
 uniformly in $k$ as $b\to\infty$.
 As a result, we conclude that (A1) holds. Q.e.d.
 
\appendix

\section{Appendix}\label{app}

{\bf Proof of Proposition \ref{pz}.} We need to show that
\begin{equation}\label{inequality}
\p(|\y_p^\top A_p \y_p-\tr(\Sigma_pA_p)|>\varepsilon p)\leqslant
\e \min\big\{16M^2\tr(\Sigma_p^2)(\varepsilon p)^{-2},1\big\}
\end{equation}
for any $\varepsilon,M>0$ and each complex $p\times p$ matrix $A_p$ with $\|A_p\|\leqslant M.$  In particular, this inequality implies that 
 (A1) holds when (A2) holds.

By Chebyshev's inequality,
\begin{align*}
\p(|\y_p^\top A_p \y_p-\tr(\Sigma_pA_p)|>\varepsilon p)=& 
\e \p(|\y_p^\top A_p \y_p-\tr(\Sigma_pA_p)|>\varepsilon p|\Sigma_p)\\
\leqslant& 
\e \min\big\{\e(|\y_p^\top A_p \y_p-\tr(\Sigma_pA_p)|^2|\Sigma_p)(\varepsilon p)^{-2},1\big\}.
\end{align*}
Write $A_p=A_{p1}+iA_{p2}$ for real $p\times p$  matrices $A_{pj},$ $j=1,2$. Then
\[|\y_p^\top A_p \y_p-\tr(\Sigma_pA_p)|^2=\sum_{j=1}^2|\y_p^\top A_{pj} \y_p-\tr(\Sigma_pA_{pj})|^2.\]
We have 
\[\y_p^\top A_{pj} \y_p=\y_p^\top A_{pj}^\top \y_p=\y_p^\top C_{pj} \y_p,\]\[\tr(\Sigma_pA_{pj})=\tr(A_{pj}^\top\Sigma_p^\top)=\tr(A_{pj}^\top\Sigma_p)=\tr(\Sigma_pA_{pj}^\top)=
\tr(\Sigma_pC_{pj}),\]
where $C_{pj}=(A_{pj}^\top+A_{pj})/2$ and $j=1,2.$ Hence,
\[|\y_p^\top A_{pj} \y_p-\tr(\Sigma_pA_{pj})|=|\y_p^\top C_{pj} \y_p-\tr(\Sigma_pC_{pj})|.\]

Let us show that $\|C_{pj}\|\leqslant M$ for $j=1,2$. Using standard properties of the spectral norm, we derive
\[\|A_{pj}\|=\sqrt{\|A_{pj}^\top A_{pj}\|}=\sqrt{\|A_{pj} A_{pj}^\top\|}=\|A_{pj}^\top\|\quad\text{and}\quad\|C_{pj}\|\leqslant\frac{\|A_{pj}\|+\|A_{pj}^\top\|}{2}=\|A_{pj}\|.\]
Thus we will prove that $\|C_{pj}\|\leqslant M$ if we show that $\|A_{pj}\|\leqslant M$ (for each $j=1,2$).

 The spectral norm $\|A_{pj}\|$ is the square root of the largest eigenvalue $\lambda$ of $A_{pj}^\top A_{pj}$. It is attained by the corresponding (real) eigenvector of $A_{pj}^\top A_{pj}$, i.e. $\|A_{pj}\|=\sqrt{\lambda}=\|A_{pj} x\|$ whenever $(A_{pj}^\top A_{pj}) x=\lambda x$ for $x\in\bR^p$ having $\|x\|=1.$ The latter implies that
\begin{align}\label{spectralnorm}
\|A_{pj}\|=\sup_{x\in\bR^p:\,\|x\|=1}&\|A_{pj}x\|\leqslant\sup_{x\in\bR^p:\,\|x\|=1}\|A_{p}x\|\leqslant\nonumber\\
&\leqslant\sup_{y\in\mathbb C^p:\,\|y\|=1}\|A_{p}y\|=\|A_{p}\|\leqslant M.\end{align}

As a result, to prove \eqref{inequality} we need to verify that 
\[\e(|\y_p^\top C_{p} \y_p-\tr(\Sigma_pC_{p})|^2|\Sigma_p)\leqslant 8M^2\tr(\Sigma_p^2)\quad\text{a.s.}\]
for   any real symmetric  $p\times p$ matrix $C_p$ with $\|C_p\|\leqslant M$. Any such matrix  can be represented as $C_p=D_{p1}-D_{p2}$ for real symmetric positive semi-definite $p\times p$ matrices $D_{pj},$ $j=1,2$, with $\|D_{pj}\|\leqslant \|C_p\|.$ Additionally, by the Cauchy inequality,
\[|\y_p^\top C_{p} \y_p-\tr(\Sigma_pC_{p})|^2\leqslant2\sum_{j=1}^2|\y_p^\top D_{pj} \y_p-\tr(\Sigma_pD_{pj})|^2.\]
Thus, to prove \eqref{inequality} we need to verify that 
\[\e(|\y_p^\top D_{p} \y_p-\tr(\Sigma_pD_{p})|^2|\Sigma_p)\leqslant 2M^2\tr(\Sigma_p^2)\quad\text{a.s.}\]
for   any real symmetric positive semi-definite $p\times p$ matrix $D_p$ with $\|D_p\|\leqslant M$. 

Let $B_p=\Sigma_p^{1/2}D_p \Sigma_p^{1/2}$. Recalling that $\y_p=\Sigma_p^{1/2}\w_p$, we have $\y_p^\top D_p \y_p=\w_p^\top B_p\w_p,$
\[\tr(\Sigma_p D_p)= \tr(B_p)\quad\text{and}\quad\tr(B_p^2)=\tr(\Sigma_pD_p \Sigma_p D_p)\leqslant \|D_p\|^2\tr(\Sigma_p^2)\]
(see \eqref{QQ}). Writing $B_p=\sum_{k=1}^p \lambda_{kp}e_{kp}e_{kp}^\top$ for orthonormal vectors $e_{kp}\in \bR^p$, $1\leqslant k\leqslant n$,  we get 
\begin{align*}
\e(|\y_p^\top &D_p \y_p-\tr(\Sigma_pD_p)|^2|\Sigma_p)=\e(|\w_p^\top B_p\w_p-\tr(B_p)|^2|\Sigma_p)=\\
=&\var\Big(\sum_{k=1}^p\lambda_{kp}(\w_p^\top e_{kp})^2|\Sigma_p\Big)=
\sum_{k=1}^p\lambda_{kp}^2\var(\xi)=2\tr(B_p^2)\leqslant 2M^2\tr(\Sigma_p^2)\end{align*}
a.s., since  $\{(\w_p^\top e_{kp})^2\}_{k=1}^p$ are independent random variables distributed as $\xi\sim \chi_1^2$ (conditionally on $\Sigma_p$).

Assume now that (A1) holds. Let us show that (A2) holds.  Take $A_p=I_p$. Hence
\[\frac{\y_p^\top \y_p-\tr(\Sigma_p)}{p}=\frac{\w_p^\top\Sigma_p \w_p-\tr(\Sigma_p)}{p}\pto 0\] and, in the above definition, $\lambda_{1p},\ldots,\lambda_{pp}$ are eigenvalues of $\Sigma_p.$ Suppose also \[\|\Sigma_p\|=\lambda_{1p}\geqslant\ldots\geqslant \lambda_{pp}\geqslant0.\] If $\w^*_p$ is an independent copy of $\w_p$ and $\w_p^*$ is also independent of $\Sigma_p$, then
\[\frac{\w_p^\top\Sigma_p \w_p-(\w_p^*)^\top\Sigma_p \w_p^*}{p}=\frac{\w_p^\top\Sigma_p \w_p-\tr(\Sigma_p)}{p}-\frac{(\w_p^*)^\top\Sigma_p \w_p^*-\tr(\Sigma_p)}{p}\pto 0.\]
Therefore,
\begin{align*}\e\exp\{i(\w_p^\top\Sigma_p \w_p-&(\w_p^*)^\top\Sigma_p \w_p^*)/p\}=\\
=&\e\prod_{k=1}^p\exp\{i \lambda_{pk}[(\w_p^\top e_{kp})^2-((\w_p^*)^\top e_{kp})^2]/p\}\\
=&\e\prod_{k=1}^p|\varphi(\lambda_{pk}/p)|^2\to 1
\end{align*}
as $p\to\infty$, where  $\varphi(t)=\e\exp\{i t\xi\}$ for $t\in\bR$ and $\xi\sim \chi_1^2$ as before. Hence, \[\e|\varphi(\|\Sigma_p\|/p)|^2=\e\frac{1}{|1-2i\|\Sigma_p\|/p|}=\e
\frac{1}{(1+4\|\Sigma_p\|^2/p^2)^{1/2}}\pto 1.\] As a result, we conclude that $\|\Sigma_p\|/p\pto0 $ and 
\begin{align*}\e\prod_{k=1}^p|\varphi(\lambda_{pk}/p)|^2
=&\e\prod_{k=1}^p\frac1{(1+4\lambda_{pk}^2/p^2)^{1/2}}\\
=&\e\min\Big\{\exp\Big\{(-2+\zeta_p)\sum_{k=1}^p\lambda_{pk}^2/p^2\Big\},1\Big\}\to 1\end{align*} for some $\zeta_p$ with $\zeta_p\pto0. $
Thus, $\tr(\Sigma_p^2)/p^2\pto0,$ i.e. (A2) holds.  Q.e.d.\\
\\
{\bf Proof of Lemma \ref{l1}.} Write $C=\sum_{k=1}^p\lambda_ke_ke_k^\top$ for orthonormal vectors $e_k\in\bR^p,$ $1\leqslant k\leqslant p.$ Then the result follows from inequalities 
\begin{align*} |1+&w^\top (C-zI_p)^{-1} w|\geqslant\Im(w^\top (C-zI_p)^{-1} w)=
\Im\Big(\sum_{k=1}^p \frac{(w^\top e_k)^2}{\lambda_k-z}\Big)=\\
&=\Im (z) \sum_{k=1}^p\frac{(w^\top e_k)^2}{|\lambda_k-z|^2}\geqslant
\Im (z) \Big|\sum_{k=1}^p\frac{(w^\top e_k)^2}{(\lambda_k-z)^2}\Big|=\Im (z)|w^\top (C-zI_p)^{-2} w|.
\end{align*}
Q.e.d.\\
{\bf Proof of Lemma \ref{l0}.} 
The spectral norm of $A$ is the square root of the largest eigenvalue of $A^*A,$ where $A^*=\overline{A^\top}=(C-\overline{z}I_p)^{-1}$. Write $z=u+iv$ for $u\in \bR$ and $v=\Im(z)>0.$ By definition,
\[A^*A=(C-\overline{z}I_p)^{-1}(C-zI_p)^{-1}=((C-uI_p)^2+v^2 I_p)^{-1}.\]
Hence, the largest eigenvalue of $A^*A$ does not exceed $1/v^2$ and $\|A\|\leqslant 1/v.$
 Q.e.d.
\\
{\bf Proof of Lemma \ref{l2}.} For any given  $\varepsilon,M>0$, set
\begin{equation}\label{IM}
I_0(\varepsilon,M)=\varlimsup_{p\to\infty}\sup_{A_p}\p(|\y_p^\top A_p \y_p-\tr(\Sigma_p A_p)|>\varepsilon p),
\end{equation} 
where the supremum is taken over all  real  symmetric $p\times p$ matrices $A_p$ with  $\|A_p\|\leqslant M.$  By this definition, there are $p_k\to\infty$ and $A_{p_k}$ with $\|A_{p_k}\|\leqslant M$  such that 
\[I_0(\varepsilon,M)=\lim_{k\to\infty}\p(|\y_{p_k}^\top A_{p_k} \y_{p_k}-\tr(\Sigma_{p_k}A_{p_k})|>\varepsilon p_k).\] 
Every real symmetric matrix $A_p$ can be represented as $A_p=A_{p1}-A_{p2}$ for real symmetric positive semi-definite $p\times p$ matrices $A_{pj},$ $j=1,2$, with $\|A_{pj}\|\leqslant \|A_p\|$. Moreover,
for any $\varepsilon>0$ and $p\geqslant 1$,
\begin{align} \label{ineq}
\p(|\y_p^\top A_p \y_p-\tr(\Sigma_p A_p)|>\varepsilon p)\leqslant & \p(|\y_p^\top A_{p1} \y_p-\tr(\Sigma_p A_{p1})|>\varepsilon p/2) \nonumber\\&+
 \p(|\y_p^\top A_{p2} \y_p-\tr(\Sigma_p A_{p2})|>\varepsilon p/2).\end{align}
Hence, it follows from (A1) that $I_0(\varepsilon,M)=0$ for any $\varepsilon,M>0.$

If $A_p$ is any real $p\times p$ matrix and $B_p=(A_p^\top +A_p)/2$, then $\y_p^\top A_p\y_p= \y_p^\top B_p\y_p$, \[
\|A_p\|=\sqrt{\|A_p^\top A_p\|}=\sqrt{\|A_p A_p^\top\|}=\|A_p^\top\|\quad\text{and}\quad\|B_p\|\leqslant\frac{\|A_p\|+\|A_p^\top\|}{2}=\|A_p\|.\]
In addition, $\tr(\Sigma_pA_p)=\tr(A_p\Sigma_p)=\tr((A_p\Sigma_p)^\top)=\tr(\Sigma_pA_p^\top )=\tr(\Sigma_pB_p).$
Thus, if $I_1(\varepsilon,M)$ is defined as $I_0(\varepsilon,M)$ in \eqref{IM} with the supremum taken over all  real $p\times p$ matrices $A_p$ with  $\|A_p\|\leqslant M$, then  
\begin{equation}\label{I1M}
I_1(\varepsilon,M)=0\quad\text{ for any $\varepsilon,M>0.$}
\end{equation}

Let now $A_p=A_{p1}+iA_{p2}$ for real $p\times p$ matrices $A_{pj},$ $j=1,2.$ It is shown in the proof of Proposition \ref{pz} (see \eqref{spectralnorm}) that 
\begin{equation}
\label{I2M}
\|A_{pj}\|\leqslant \|A_p\|,\quad j=1,2.\end{equation} Define $I_2(\varepsilon,M)$ similarly to  $I_0(\varepsilon,M)$ in \eqref{IM} with the supremum taken over all  complex $p\times p$ matrices $A_p$ with  $\|A_p\|\leqslant M$. Combining \eqref{ineq}, \eqref{I1M} and \eqref{I2M} yields $I_2(\varepsilon,M)=0$ for any $\varepsilon,M>0.$ Q.e.d.
\\
{\bf Proof of Lemma \ref{l3}.} We have
\begin{align*} I=\frac{z_1}{1+w_1}-\frac{z_2}{1+w_2}&=\frac{z_1(1+w_2)-z_2(1+w_1)-z_1w_1+w_1z_1}{(1+w_1)(1+w_2)}\\&=
\frac{(z_1-z_2)+z_1(w_2-w_1)+w_1(z_1-z_2)}{(1+w_1)(1+w_2)}.\end{align*}
It follows from $|z_1-z_2|\leqslant \gamma$, $|w_1-w_2|\leqslant \gamma$ and $|z_1|/|1+w_1|\leqslant M$ that 
\[|I|\leqslant \frac{\gamma(1+|z_1|+|w_1|)}{|1+w_1||1+w_2|}\leqslant \frac{\gamma}{|1+w_2||1+w_1|}+\frac{\gamma M}{|1+w_2|}+\frac{\gamma}{|1+w_2|}\,\frac{|w_1|}{|1+w_1|}.\]
In addition, we have  $|1+w_2| \geqslant \delta,$
\[|1+w_1|=|1+w_2+(w_1-w_2)|\geqslant \delta-\gamma\geqslant \delta /2,\]
\[
\frac{|w_1|}{|1+w_1|}= \dfrac{2}{|1+w_1|} I(|w_1|\leqslant 2)+
\dfrac{|w_1|}{|w_1|-1}\,I(|w_1|>2)\leqslant\begin{cases}
4/\delta,& |w_1|\leqslant 2,\\
2,& |w_1|>2.\end{cases}\]
Finally, we conclude that  $|I|\leqslant \gamma(M/\delta +4/\min\{\delta^2,2\delta\}+2/\delta^2).$ Q.e.d.
\\
{\bf Proof of Lemma \ref{l4}.} Write $z=u+iv$ for $u\in \bR$ and $v>0$. We need to prove inequality
\[|z||1+\tr(\Sigma(C-zI_p)^{-1})|=|z+\tr(\Sigma(C/z-I_p)^{-1})|\geqslant v.\]
Since 
\[|z+\tr(\Sigma(C/z-I_p)^{-1})|\geqslant \Im(z+\tr(\Sigma(C/z-I_p)^{-1}))=v+\Im(\tr(\Sigma(C/z-I_p)^{-1})),\]
we only need to check that
\[\Im\big(\tr(\Sigma(C/z-I_p)^{-1})\big)\geqslant0.\]
Denote
\begin{equation} \label{B}B=\Big(\frac{u}{|z|^2}C-I_p\Big)^2+\frac{v^2}{|z|^4}C^2.\end{equation}
Note that $B$ is an invertible symmetric positive definite matrix, since \[B=(C/z-I_p)(C/z-I_p)^*=(C/z-I_p)^*(C/z-I_p)=\]\[=
(C/\overline{z}-I_p)(C/z-I_p)=\frac{1}{|z|^2} C^2-\frac{2u}{|z|^2} C+I_p\] and $C/z-I_p=(C-zI_p)/z$ is invertible, where $A^*=\overline{A^\top}$ is the conjugate transpose of a matrix $A$.
Additionally, 
\begin{align*}
(C/z-I_p)^{-1}=(C/\overline{z}-I_p)B^{-1}=\Big(\frac{u}{|z|^2}C-I_p+\frac{iv}{|z|^2} C\Big)B^{-1}.
\end{align*}
Therefore,
\[\Im\big(\tr(\Sigma(C/z-I_p)^{-1})\big)=\frac{v}{|z|^2}\Im\big(\tr(\Sigma CB^{-1})\big).\]
Let $C^{1/2}$ and $\Sigma^{1/2}$ be the principal square roots of $C$ and $\Sigma$. Then, by the definition of $B$, matrices $C^{1/2}$ and $B^{-1}$ commute, $CB^{-1}=C^{1/2}B^{-1}C^{1/2}$ and
\[\tr(\Sigma CB^{-1})=\tr(\Sigma C^{1/2}B^{-1}C^{1/2})=\tr(\Sigma^{1/2} C^{1/2}B^{-1}C^{1/2}\Sigma^{1/2}).\]
As it is shown above, $B$ is symmetric positive definite. Hence, $B^{-1}$ is symmetric positive definite and 
$QB^{-1}Q^\top$ is symmetric positive semi-definite for any $p\times p$ matrix $Q$. Taking $Q=\Sigma^{1/2} C^{1/2}$, we see that 
 \[\tr(\Sigma CB^{-1})=\tr(QB^{-1}Q^\top)\geqslant 0.\]
 This proves the desired bound. Q.e.d.
\\
{\bf Proof of Lemma \ref{l8}.} It is shown in the proof of Lemma \ref{l0} that, for  all $\varepsilon,M>0$ and  each $k=1,\ldots,n,$
\begin{align*}
\sup_{A_p}\p(|\y_{pk}^\top A_p \y_{pk}-\tr(\Sigma_{pk} A_p)|>\varepsilon p)&\leqslant 
2\sup_{B_p}\p(|\y_{pk}^\top B_p \y_{pk}-\tr(\Sigma_{pk} B_p)|>\varepsilon p/2)\\
&\leqslant 
2\sup_{C_p}\p(|\y_{pk}^\top C_p \y_{pk}-\tr(\Sigma_{pk} C_p)|>\varepsilon p/2)\\&\leqslant 
4\sup_{D_p}\p(|\y_{pk}^\top D_p \y_{pk}-\tr(\Sigma_{pk} D_p)|>\varepsilon p/4),
\end{align*}
where $A_p$ is a complex $p\times p$ matrix with $\|A_p\|\leqslant M$,  $B_p$ is a real $p\times p$ matrix with $\|B_p\|\leqslant M$, $C_p$ is a real symmetric $p\times p$ matrix with $\|C_p\|\leqslant M,$ and $D_p$ is a real symmetric positive semi-definite $p\times p$ matrix with $\|D_p\|\leqslant M.$ Thus the result follows from (A3). Q.e.d.
\\
{\bf Proof of Lemma \ref{l5}.} In what follows, we denote the principal square root of $Q$ by $Q^{1/2}$ for  any real symmetric positive semi-definite $p\times p$ matrix $Q$. For each $z=u+iv$ with $u\in\bR$ and $v=\Im(z)>0$,
\[\|(I_q+U^\top(C-zI_p)^{-1}U)^{-1}\|=|z|\|(zI_q+U^\top(C/z-I_p)^{-1}U)^{-1}\|.\]
Arguing as in the proof of Lemma \ref{l4}, we derive
\[(C/z-I_p)^{-1}=\Big(\frac{u}{|z|^2}C-I_p\Big) B^{-1}+\frac{iv}{|z|^2} CB^{-1}=A_1+iA_2,\]
where $B$ is given in \eqref{B}, $A_1$ and $A_2$ are real symmetric commuting $p\times p$ matrices defined by
\[A_1=\frac{u}{|z|^2}C^{1/2}B^{-1}C^{1/2}-B^{-1}, \quad A_2=\frac{v}{|z|^2} C^{1/2}B^{-1}C^{1/2}.\]
In addition, $A_2$ is symmetric positive semi-definite. We have
\begin{align*}\|(zI_q+U^\top(C/z-I_p)^{-1}U)^{-1}\|
&=\|(  (uI_q+U^\top A_1 U)+i(vI_q+U^\top A_2 U))^{-1}\|\\
&\leqslant \big\|A_3^{-1/2}(A_3^{-1/2}(uI_q+U^\top A_1 U)A_3^{-1/2}+iI_q)^{-1}A_3^{-1/2}\big\|\\
&\leqslant \big\|A_3^{-1/2}\big\|^2\big\|(A_3^{-1/2}(uI_q+U^\top A_1 U)A_3^{-1/2}+iI_q)^{-1}\big\|,
\end{align*}
where $A_3=vI_p+U^\top A_2 U$. Since $U^\top A_2 U$ is a symmetric positive semi-definite matrix,
the spectral norm of $A_3^{-1/2}$ does not exceed $1/\sqrt{v}.$ The spectral norm of 
\[A_4=(A_3^{-1/2}(uI_q+U^\top A_1 U)A_3^{-1/2}+iI_q)^{-1}\]
is the largest eigenvalue of 
\[A_4^*A_4=\big(\big(A_3^{-1/2}(uI_q+U^\top A_1 U)A_3^{-1/2}\big)^2+I_q\big)^{-1}.\]
Obviously, it is not greater than $1$. Finally, we conclude that
\[\|(I_q+U^\top(C-zI_p)^{-1}U)^{-1}\|\leqslant \frac{|z|}v.\]
 Q.e.d.
\\
{\bf Proof of Lemma \ref{l6}.} As above, we denote the principal square root of $Q$ by $Q^{1/2}$ for  any real symmetric positive semi-definite $p\times p$ matrix $Q$.  Let also $R^*=\overline{R^\top}$ be  the conjugate transpose of a  matrix (or a vector) $R$.

Write $z=u+iv$ for $u\in\bR$ and $v>0$ and set
\[A_0=(C-zI_p)^{-1}U(I_q+U^\top(C-zI_p)^{-1}U)^{-1}.\]
We need to prove that 
\[|w^\top A_0^\top A_0 w|\leqslant \frac{v+|z|}{v^2} \|w\|^2.\]
By the Cauchy-Schwartz inequality,
\[|w^\top A_0^\top A_0 w|=\big|\big(\overline{A_0w},A_0w\big)\big|\leqslant\big\|\overline{A_0w}\big\|\,\|A_0w\|=\|A_0w\|^2\leqslant \|A_0\|^2 \|w\|^2.\]

By Theorem A.6 in \cite{BS},
\begin{align*}
(C-zI_p)^{-1}=(C-uI_p&-iv I_q)^{-1}=(C-uI_p+iv I_p)B^{-1}=\\
&=CB^{-1}-\overline{z}B^{-1}=C^{1/2}B^{-1}C^{1/2}-\overline{z}B^{-1},
\end{align*}
\[[(C-zI_p)^{-1}]^*(C-zI_p)^{-1}=(C-\overline{z}I_p)^{-1}(C-zI_p)^{-1} =B^{-1},\]
where \begin{align}\label{hz}
B=(C-uI_p)^2+v^2I_p\end{align} is a real symmetric positive definite matrix  commuting with $C^{1/2}$ (its inverse $B^{-1}$ has the same properties). We also have
\begin{align*}
(I_q+U^\top(C-zI_p)^{-1}U)^{-1}&=(A_1-\overline{z}\,U^\top B^{-1}U)^{-1}\\
&=A_1^{-1/2}(I_q-\overline{z}\,A_1^{-1/2}U^\top B^{-1}UA_1^{-1/2})^{-1}A_1^{-1/2}\\
&=A_1^{-1/2}(I_q-\overline{z}V^\top V)^{-1} A_1^{-1/2}
\end{align*}
and  $A_0=A_2A_1^{-1/2},$
where  $V=B^{-1/2}U A_1^{-1/2}$, \begin{align}\label{A1}A_1&=I_q+U^\top C^{1/2}B^{-1}C^{1/2} U,\\
A_2&=(C-zI_p)^{-1}B^{1/2} V (I_p-\overline{z} V^\top V)^{-1}.\nonumber
\end{align}  The matrix $Q^\top B^{-1}Q$ is real symmetric positive semi-definite for any $p\times p$ matrix $Q$. Therefore, taking $Q=C^{1/2} U$ yields \[\|A_1^{-1/2}\|=\|(I_q+Q^\top B^{-1} Q)^{-1/2}\|=\|(I_q+Q^\top B^{-1} Q)^{-1}\|^{1/2}\leqslant 1\] and $\|A_0\|\leqslant \|A_2\|\,\|A_1^{-1/2}\|\leqslant 
\|A_2\|.$ The spectral norm $\|A_2\|$ is equal to the square root of the spectral norm of the matrix
\begin{align*}
 A_2^*A_2&=[(I_p-\overline{z} V^\top V)^{-1}]^*V^\top B^{1/2} [(C-zI_p)^{-1}]^*(C-zI_p)^{-1}B^{1/2}V (I_p-\overline{z} V^\top V)^{-1}\\
 &=(I_p-z V^\top V)^{-1}V^\top V (I_p-\overline{z} V^\top V)^{-1}.
 \end{align*}
 Additionaly, 
\begin{align*}
\| A_2^*A_2\|&\leqslant \frac{1}{|z|^2}\|(V^\top V-\overline{z}I_p/|z|^2)^{-1}V^\top V (V^\top V-zI_p/|z|^2)^{-1}\|\\
&\leqslant \frac{1}{|z|^2}\|(V^\top V-\overline{z}I_p/|z|^2)^{-1}\| \| V^\top V (V^\top V-zI_p/|z|^2)^{-1}\|.
 \end{align*}
By the triangle inequality,
\begin{align*} \| V^\top V (V^\top V-zI_p/|z|^2)^{-1}\|=&\|I_p+z|z|^{-2}(V^\top V-zI_p/|z|^2)^{-1}\|\\&\leqslant 1+
\frac{1}{|z|}\|(V^\top V-zI_p/|z|^2)^{-1}\|
\end{align*}
Matrices $A_3=(V^\top V-zI_p/|z|^2)^{-1}$ and $A_3^*=(V^\top V-\overline{z}I_p/|z|^2)^{-1}$ have identical spectral norms equal to
the square root of the largest eigenvalue of the matrix 
\[A_3^*A_3=A_3A_3^*=((V^\top V-uI_p/|z|^2)^{2}+v^2 I_p/|z|^4)^{-1}.\]
Obviously, this eigenvalue does not exceed $|z|^4/v^2.$ Combining the above estimates yields
\[|w^\top A_0^\top A_0 w|\leqslant\|A_0\|^2\|w\|^2\leqslant \frac{1}{|z|^2} \,\frac{|z|^2}{v}\Big(1+\frac1{|z|}\,\frac{|z|^2}{v}\Big)\|w\|^2=\frac{v+|z|}{v^2} \|w\|^2.\]

\noindent{\bf Proof of Lemma \ref{l7}.} Let $z=u+iv$ for $u\in\bR$, $v=\Im(z)>0$ and $C_z=C-zI_p$. It follows from Lemma \ref{l5} that \[A_z=I_q+U^\top C_z^{-1} U\] is non-degenerate. By  the Sherman-Morrison-Woodbury formula,
\begin{equation}\label{smw}
(C_z+UU^\top)^{-1}= C_z^{-1}-C_z^{-1}UA_z^{-1} U^\top C_z^{-1}.\end{equation}
The Cauchy-Schwartz inequality, \eqref{smw} and Lemma \ref{l5} imply that 
\begin{align*}
|y^\top(C_z+UU^\top)^{-1}y-y^\top C_z^{-1}y|&\leqslant (y^\top C_z^{-1}U)A_z^{-1}(U^\top C_z^{-1}y)\\
&\leqslant \|\overline{U^\top C_z^{-1}y}\|\,\|A_z^{-1}(U^\top C_z^{-1}y)\|\\
&\leqslant \frac{|z|}v\|U^\top C_z^{-1}y\|^2.
\end{align*}
It also follows from \eqref{smw} that
\begin{align*} 
y^\top(C_z+UU^\top)^{-2}y=&y^\top( C_z^{-1}-C_z^{-1}UA_z^{-1}U^\top C_z^{-1})^2y\\
=&y^\top C_z^{-2}y+ (y^\top C_z^{-1}U)A_z^{-1}U^\top C_z^{-2}UA_z^{-1}(U^\top C_z^{-1}y)-\\
&-(y^\top C_z^{-1}U)A_z^{-1}(U^\top C_z^{-2}y)-(y^\top C_z^{-2}U)A_z^{-1}(U^\top C_z^{-1}y).
\end{align*}
Applying the Cauchy-Schwartz inequality, this identity, Lemma \ref{l5} and Lemma \ref{l6}, we get 
\begin{align*}
|y^\top(C_z+UU^\top)^{-2}y-&y^\top C_z^{-2}y|\leqslant 
\frac{2|z|}v\|U^\top C_z^{-1}y\|\,\|U^\top C_z^{-2}y\| +\frac{v+|z|}{v^2} \|U^\top C_z^{-1}y\|^2\\
&\leqslant \frac{|z|}v\|U^\top C_z^{-1}y\|^2+\frac{|z|}v\|U^\top C_z^{-2}y\|^2 +\frac{v+|z|}{v^2} \|U^\top C_z^{-1}y\|^2\\
&\leqslant \frac{(|z|+1)^2}{v^2}\sum_{j=1}^2\|U^\top C_z^{-j}y\|^2.
\end{align*}
Gathering together the above estimates, we finish the proof. Q.e.d.

\noindent{\bf Proof of Lemma \ref{l10}.}  Let $z=u+iv$ for $u\in\bR$ and $v=\Im(z)>0.$ 
 Denote further by $Q^*$ the conjugate transpose of a matrix $Q$, i.e. $Q^*=\overline{Q^\top}$. 

 Recall that $\tr(Q_1Q_2)=\tr(Q_2Q_1),$ $(Q_1Q_2)^*=Q_2^*Q_1^*$, $(Q^{-1})^*=(Q^*)^{-1}$, $(Q^*)^*=Q,$
 \[\|Q^*\|=\sqrt{\|QQ^*\|}=\sqrt{\|Q^*Q\|}=\|Q\|,\] \begin{equation}\label{Dnorm}
 |\tr(Q)|^2= |\tr(I_q^*Q)|^2\leqslant \tr(I_q^*I_q)\tr(Q^*Q)= q\,\tr(Q^*Q)\leqslant q^2\|Q\|^2
 \end{equation}  for any complex $q\times q$ matrices $Q,Q_1,Q_2$.  Hence, taking $Q=AU(I_q+U^\top AU)^{-1}U^\top A$ we arrive at the bound
 \[|\tr(U^\top A^2U(I_q+U^\top AU)^{-1})|^2=|\tr(Q)|^2\leqslant q\,\tr(Q^*Q),\]
 where
 \begin{align*}\tr&(Q^*Q)=\tr(A^* U(I_q+U^\top A^*U)^{-1}U^\top A^*A U(I_q+U^\top AU)^{-1}U^\top A)\\
 &=\tr(A_1^{-1/2}U^\top AA^* U(I_q+U^\top A^*U)^{-1}A_1^{1/2}A_1^{-1/2}U^\top A^*A U(I_q+U^\top AU)^{-1} A_1^{1/2})
 \end{align*}
 and $A_1$ is given in \eqref{A1}. Thus, $\tr(Q^*Q)$ is bounded from above by 
\begin{align*} 
I=&q\|A_1^{-1/2}U^\top AA^* U(I_q+U^\top A^*U)^{-1}A_1^{1/2}\|\|A_1^{-1/2}U^\top A^*A U(I_q+U^\top AU)^{-1} A_1^{1/2}\|\\
=&q\|A_1^{1/2}(I_q+U^\top AU)^{-1} U^\top AA^* U A_1^{-1/2}\|\|A_1^{-1/2}U^\top A^*A U(I_q+U^\top AU)^{-1} A_1^{1/2}\|
\end{align*}
 Since $A=(C-zI_p)^{-1}$, $A^*=(C-\overline{z}I_p)^{-1}$ and $A^*A=AA^*=B^{-1}$ for $B$ defined in \eqref{hz}, we can proceed further as in the proof of Lemma \ref{l6} using the same notation. Namely,
 \[A_1^{-1/2}U^\top A^*A U(I_q+U^\top AU)^{-1}A_1^{1/2}=V^\top V(I_q-\overline{z}V^\top V)^{-1},\]
 \[A_1^{1/2}(I_q+U^\top A U)^{-1}U^\top AA^* UA_1^{-1/2}=(I_q-\overline{z}V^\top V)^{-1}V^\top V.\]
 By the same arguments as in the proof of Lemma \ref{l6},
 \[\|V^\top V(I_q-\overline{z}V^\top V)^{-1}\|\leqslant \frac{1}{|z|}\Big(1+\frac{|z|}{v}\Big),\]
 \[\|(I_q-\overline{z}V^\top V)^{-1}V^\top V\|\leqslant \frac{1}{|z|}\Big(1+\frac{|z|}{v}\Big).\]
Finally we conclude that
\[|\tr(U^\top A^2U(I_q+U^\top AU)^{-1})|\leqslant \frac{q(|z|+v)}{|z|v}.\]
Q.e.d.

\noindent{\bf Proof of Lemma \ref{l9}.} Let $z=u+iv$ for $u\in\bR$ and $v=\Im(z)>0,$ 
\[\Delta=V^\top AA^*V-U^\top AA^*U\quad\text{and}\quad\Delta_j=V^\top A^{j}V-U^\top A^{j}U,\quad j=1,2.\] Denote further by $Q^*$ the conjugate transpose of a matrix $Q$, i.e. $Q^*=\overline{Q^\top}$.  

We have
\begin{align*}
|\tr\big(V^\top &A^{2}V(I_q+V^\top AV)^{-1}\big)-\tr\big(U^\top A^{2}U(I_q+U^\top AU)^{-1}\big)|\leqslant\\
\leqslant&|\tr\big(V^\top A^2 V(I_q+V^\top AV)^{-1}\big)-\tr\big(U^\top A^{2}U(I_q+V^\top AV)^{-1}\big)|\\
&+|\tr\big(U^\top A^{2}U(I_q+V^\top AV)^{-1}\big)-\tr\big(U^\top A^{2}U(I_q+U^\top AU)^{-1}\big)|\\
\leqslant & |\tr\big(\Delta_2(I_q+V^\top AV)^{-1}\big)|+|\tr\big(U^\top A^{2}U(I_q+U^\top AU)^{-1}\Delta_1(I_q+V^\top AV)^{-1}\big)|
 \end{align*}
 By Lemma \ref{l5} and \eqref{Dnorm}, \begin{align*}
 |\tr\big(\Delta_2(I_q+V^\top AV)^{-1}\big)|\leqslant& q\|\Delta_2(I_q+V^\top AV)^{-1}\|\\
 \leqslant& q\|\Delta_2\|\|(I_q+V^\top AV)^{-1}\|\\
 \leqslant& \frac{q|z|}{v}\|\Delta_2\|.
 \end{align*}
 Taking $Q=AU(I_q+U^\top AU)^{-1}\Delta_1(I_q+V^\top AV)^{-1}U^\top A$ we infer that
 \[|\tr\big(U^\top A^{2}U(I_q+U^\top AU)^{-1}\big)\Delta_1(I_q+V^\top AV)^{-1}\big)|=|\tr(Q)|\leqslant q\|Q\|\]
and
\[\|Q\|\leqslant\|AU(I_q+U^\top AU)^{-1}\|\|\Delta_1\|\|(I_q+V^\top AV)^{-1}U^\top A\|.\]
It is shown in the proof of Lemma \ref{l6} (where $A_0=AU(I_q+U^\top AU)^{-1}$) that 
\begin{equation} \label{U+}
\|AU(I_q+U^\top AU)^{-1}\|\leqslant\frac{\sqrt{v+|z|}}{v}. \end{equation}
In addition,  since $\|R\|^2=\|R^*R\|=\|RR^*\|$ for any complex $q\times p$ matrix $R,$ we infer that 
\begin{align*} \|(I_q+&V^\top AV)^{-1}U^\top A\|^2=
\|(I_q+V^\top AV)^{-1}U^\top AA^*U[(I_q+V^\top AV)^{-1}]^*\|\\
\leqslant &\|(I_q+V^\top AV)^{-1}\Delta[(I_q+V^\top AV)^{-1}]^*\|\\
&+\|(I_q+V^\top AV)^{-1}V^\top AA^*V[(I_q+V^\top AV)^{-1}]^*\|\\
\leqslant&\|(I_q+V^\top AV)^{-1}\|^2\|\Delta\|+\|(I_q+V^\top AV)^{-1}V^\top A\|^2.
\end{align*}
 
As in the proof of Lemma \ref{l6}, one can show that
\[\|(I_q+V^\top AV)^{-1}V^\top A\|\leqslant \frac{\sqrt{v+|z|}}{v}.\] 
By Lemma \ref{l5},
\[\|(I_q+V^\top AV)^{-1}\|\leqslant \frac{|z|}{v}.\]
Combining the above estimates yields
\begin{align*}
|\tr(Q)|\leqslant &q\frac{\sqrt{v+|z|}}{v}\|\Delta_1\|\bigg(\frac{|z|^2}{v^2}\|\Delta\|+
\frac{v+|z|}{v^2}\bigg)^{1/2}\\
\leqslant &\frac{q|z|\sqrt{v+|z|}}{v^2}\|\Delta_1\| \sqrt{\|\Delta\|}+
\frac{q(v+|z|)}{v^2}\|\Delta_1\|\\
\leqslant &\frac{2q(|z|+1)^{3/2}}{v^2}(\|\Delta_1\| \sqrt{\|\Delta\|}+
\|\Delta_1\|).
\end{align*}
The latter gives the desired bound
\begin{align*}
|\tr\big(V^\top &A^{2}V(I_q+V^\top AV)^{-1}\big)-\tr\big(U^\top A^{2}U(I_q+U^\top AU)^{-1}\big)|\leqslant \\
&\leqslant \frac{2q(|z|+1)^{3/2}}{v^2}(\|\Delta_1\| \sqrt{\|\Delta\|}+\|\Delta_1\|+\|\Delta_2\|).
\end{align*}
Q.e.d.

\end{document}